\documentclass[10pt]{amsart}
\usepackage{amssymb,amsmath,amscd,amsfonts,amsthm,mathrsfs}
\usepackage{verbatim, array} 
\usepackage{anysize}
\usepackage{graphicx}
\usepackage{calc}
\newsavebox\CBox
\newcommand\hcancel[2][0.5pt]{%
  \ifmmode\sbox\CBox{$#2$}\else\sbox\CBox{#2}\fi%
  \makebox[0pt][l]{\usebox\CBox}%
  \rule[0.5\ht\CBox-#1/2]{\wd\CBox}{#1}}

\input xy
\xyoption{all}
\usepackage[all]{xy}  
\usepackage{color}

\usepackage[latin1]{inputenc}

\newcommand{\N}{\mathbb{N}}
\newcommand{\Z}{\mathbb{Z}}

\newcommand{\R}{\mathbb{R}}

\newcommand{\Bc}{\mathcal{B}}
\newcommand{\Cc}{\mathcal{C}}

\newcommand{\Ec}{\mathcal{E}}
\newcommand{\Gc}{\mathcal{G}}

\newcommand{\Pc}{\mathcal{P}}

\newcommand{\Tc}{\mathcal{T}}

\def\build#1_#2^#3{\mathrel{\mathop{\kern 0pt#1}\limits_{#2}^{#3}}}
\newtheorem{theorem}{Theorem}[section]
\newtheorem{proposition}[theorem]{Proposition}
\newtheorem{lemma}[theorem]{Lemma}
\newtheorem{remark}[theorem]{Remark}
\newtheorem{example}[theorem]{Example}
\newtheorem{definition}[theorem]{Definition}
\newtheorem{corollary}[theorem]{Corollary}
\newtheorem{conjecture}[theorem]{Conjecture}

\numberwithin{equation}{section}

\newcommand{\beq}{\begin{equation}}
\newcommand{\eeq}{\end{equation}}

\newcommand{\beas}{\begin{eqnarray*}}
\newcommand{\eeas}{\end{eqnarray*}}

\def \lint{[\![}
\def \rint{]\!]}
\def \biglint{\big[\!\!\big[}
\def \bigrint{\big]\!\!\big]}

\title{The closed knight tour problem in higher dimensions}

\begin{document}
\author{Joshua Erde}
\address{Department of Pure Mathematics and Mathematical Statistics, Centre for Mathematical
Sciences,
\\Wilberforce Road, Cambridge, CB3 0WB, UK}
\email{jpe28@cam.ac.uk}
\author{Bruno Gol\'enia}
\address{Department of Computer Science, University of Bristol,
  Merchant Venturers Building, Woodland Road 
\\BRISTOL BS8 1UB, United Kingdom}
\email{goleniab@compsci.bristol.ac.uk}
\author{Sylvain Gol\'enia}
\address{Institut de Math\'ematiques de Bordeaux, Universit\'e
Bordeaux $1$, $351$, cours de la Lib\'eration
\\$33405$ Talence cedex, France}
\email{sylvain.golenia@u-bordeaux1.fr}

\date{Version of \today}
\begin{abstract}
The problem of existence of closed knight tours for rectangular chessboards was solved by Schwenk in 1991. Last year, in 2011, DeMaio and Mathew provide an extension of this result for $3$-dimensional rectangular boards. 
In this article, we give the solution for $n$-dimensional rectangular boards,
for $n\geq 4$. 
\end{abstract}

\maketitle

\tableofcontents

\section{Introduction}
On a chessboard, a knight moves by two squares in one direction and by
one square in the other one (like a L). A classical challenge is the
so-called \emph{knight tour}. The knight is placed on the empty board
and, moving according to the rules of chess, must visit each square
exactly once. A knight's tour is called a \emph{closed tour} if the knight
ends on a square attacking the square from which it began. If the
latter is not satisfied and the knight has visited each square
exactly once, we call it an \emph{open tour}. 

\begin{figure}[ht]
\centering
\includegraphics[scale=0.2]{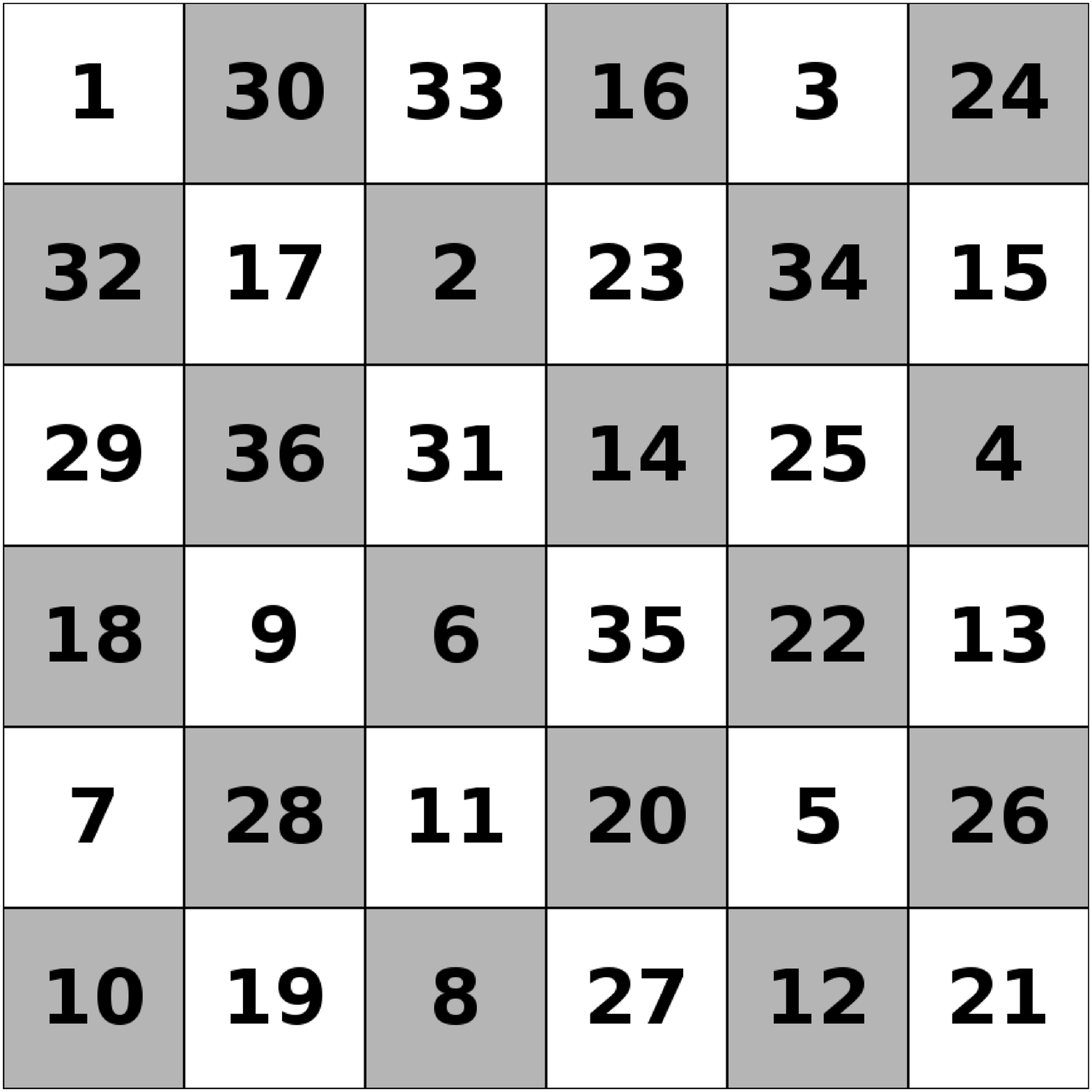}
\caption{A $6\times 6$ closed tour}\label{f:tour}
\end{figure}

Some early
solutions were given by Euler, see \cite{Eul} and also by De Moivre
(we refer to Mark R. Keen for historical remarks, see \cite{Kee}). The
problem was recently considered for various types of 
chessboards: as a cylinder \cite{Wat1},  a torus \cite{Wat2}, a
sphere \cite{Cai}, the exterior of the cube \cite{QiW}, the interior
of the cube \cite{De},$\ldots$ It represents also an active field of
research in computer science, e.g., \cite{Pab} (see references
therein). In this paper, we shall focus on rectangular boards.   

In 1991, Schwenk considered the question of the closed knight tour
problem in  a $2$-dimensional rectangular chessboard. He provided a
necessary and sufficient condition on the size of the board in order
to have a closed knight tour. He obtained: 

\begin{theorem}[Schwenk]\label{t:S} Let $1\leq m\leq n$. The $m\times n$
  chessboard has \emph{no} closed knight tour if and only if
one of the following assumption holds:
\begin{enumerate}
\item $m$ and $n$ are both odd,
\item $m\in\{1,2,4\}$,
\item $m=3$ and $n\in \{4,6,8\}$.
\end{enumerate} 
\end{theorem}

We refer to \cite{Sch} (see also
\cite{Wat}) for a proof.  When conditions (a),(b), and (c)
are not fulfilled, he reduced the problem to studying a finite
number of elementary boards. On each of them, he exhibited a
closed tour and then explained how to ``glue" the elementary boards together, in
order to make one closed tour for the union out of the disjoint
ones given by the elementary blocks. We explain the latter on an
example. Say we want a closed tour for a $12\times 6$
board. Write, side by side,  two copies  of Figure \ref{f:tour}. Delete the
connection between $21$ and $22$ for the left board and the connection
between $28$ and $29$ for the right one. Then link $21$ with $28$ and
$22$ with $29$. The Hamiltonian cycle goes as follows:  
\begin{align*}
1[L] &\to 2[L] \to \ldots \to  21[L] \to 28[R] \to 27 [R] \to \ldots \to 
1[R]\to
\\
&\to 36[R] \to \ldots \to 29 [R] \to 22[L] \to 23 [L] \to \dots \to
36[L] \to 1 [L], 
\end{align*}
where $[L]$ and $[R]$ stand for left and right, respectively.   

We turn to the question for higher dimensions. In dimension
$3$ or above, a knight moves by two steps along one coordinate
and by one step along a different one. We refer to Section
\ref{s:notation} for a mathematical definition. Stewart \cite{Ste} and
DeMaio \cite{De} constructed some examples of $3$-dimensional knight
tours. Then, in 2011, in \cite{DeM}, DeMaio and Mathew extended Theorem
\ref{t:S} by classifying all the $3$-dimensional rectangular
chessboards which admit a knight tour.  

\begin{theorem}[DeMaio and Mathew]\label{t:DeM} Let $2\leq m\leq n\leq
  p$.  The $m\times   n\times p$   chessboard has \emph{no} closed
  knight tour if and only   if one of the following assumption holds:
\begin{enumerate}
\item $m$, $n$, and $p$ are all odd,
\item $m=n=2$,
\item $m=2$ and $n=p=3$.
\end{enumerate} 
\end{theorem}
The strategy is the same as in Theorem \ref{t:S}. We give an
alternative proof of this result in Appendix \ref{s:DeMalt}.

\begin{figure}[ht]
\centering
\includegraphics[scale=0.2]{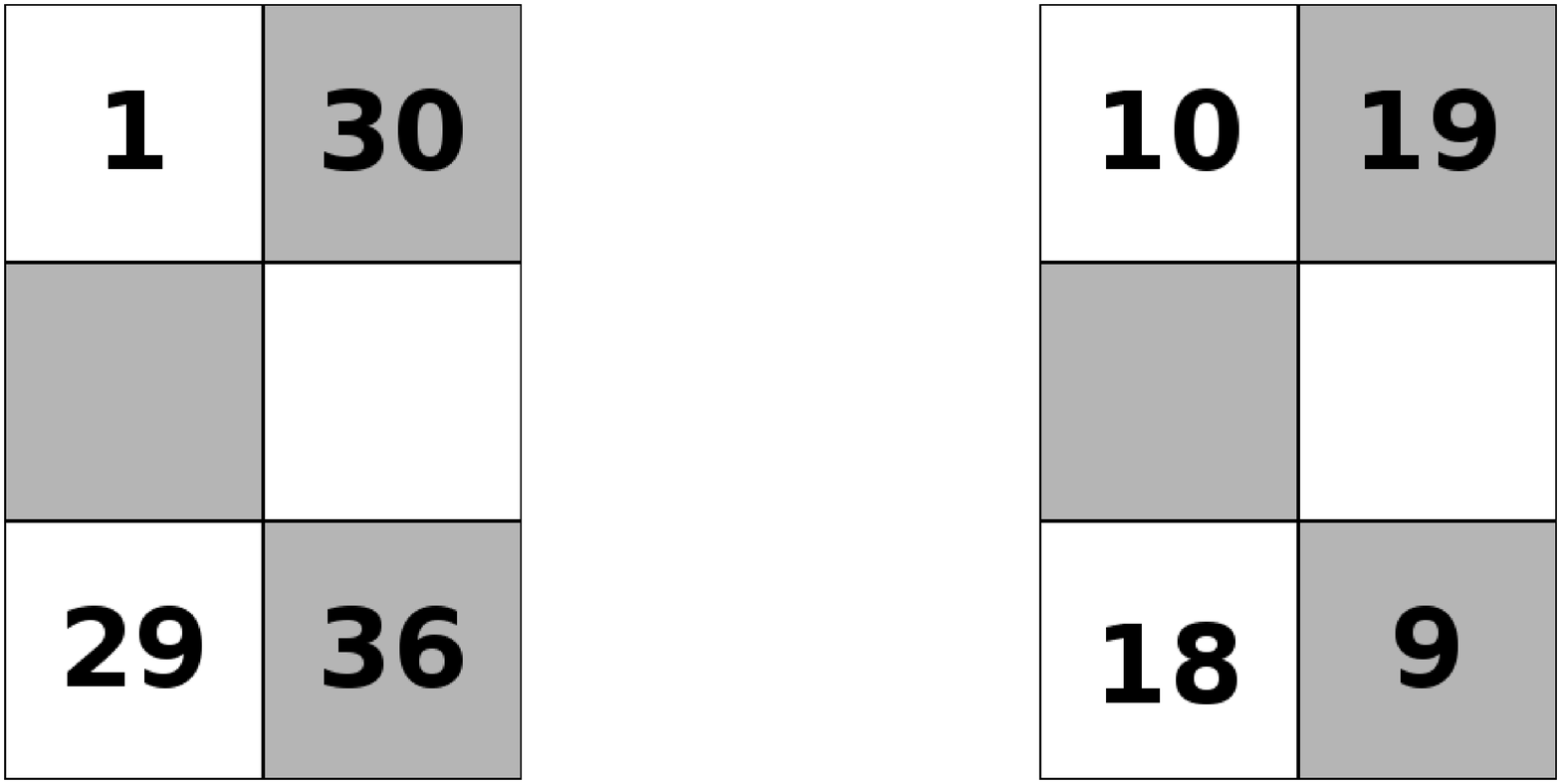}
\caption{Two cross-patterns}\label{f:2cp}
\end{figure}

In this paper we extend the previous results to higher dimensional
boards. We rely strongly on the structure of the solutions for the
case $n=3$ to treat the case $n\geq 4$. We proceed by
induction. Before giving
the main statement, we explain the key idea with Figure \ref{f:tour}. We
first extract two \emph{cross-patterns} (cp).  We
represent them up to some rotation, see Figure \ref{f:2cp}. Note they
are with disjoint support.  

We construct a tour for a $6\times 6\times 2$ board first. We work
with coordinates. The tour given in \ref{f:tour} is given by
$(a^i)_{i\in \lint 1, 36 \rint }$, with $a^1:=(1, 6)$, $a^2:=(3, 5)$,
$\ldots$ We take two copies of the tour given in
\ref{f:tour} and denote them 
by $(a^i, 1)_{i\in \lint 1, 36 \rint}$ and $(a^i, 2)_{i\in \lint 1,
  36 \rint }$ for the first and second copy, respectively. 
We can construct a tour as follows:
\begin{align*}
(3,2, 1)\to (1,2, 2) \to (3,3, 2) \to \ldots \to (3,1, 2)
\to (1,1, 1) \to \ldots \to (3,2, 1).
\end{align*}
We point out that we use coordinates in matrix way. To enhance
the idea we rewrite it abusively by 
\begin{align*}
(``36", 1)\to (``30", 2) \to (``31", 2) \to \ldots \to (``29", 2)
\to (``1", 1) \to \ldots \to (``36", 1).
\end{align*}
We explain how to treat a $6\times 6\times k$ board for $k\geq 3$. As
the proof is the same, we write the case $k=3$. We use the two cp as follows:
\begin{align}\nonumber
(``36", 1)&\to (``30", 2) \to (``31", 2) \to \ldots \to (``9", 2)\to
\\
\nonumber
&\to (``19", 3)\to (``20", 3)\to \ldots \to (``18", 3) \to 
\\
\label{e:tour2}
&\to (``10", 2)\to (``11", 2)\to \ldots \to (``29", 2) \to
\\
\nonumber
&\to (``1", 1) \to \ldots \to (``36", 1).
\end{align}
Note that we have use only one cp on the first copy and one  on last
one. Two of them are still free. We can therefore repeat the
procedure and add inductively further dimensions. For instance, we can
go from a tour for a $6\times 6\times k$ board to one for a $6\times
6\times k\times l$ board, with $k, l\geq 2$. One takes $l$ copies
of the tour. By noticing that there are two cp on the initial tour,
we proceed as in \eqref{e:tour2} for the tour. Thus, we will get a
tour for the $6\times 6\times k\times l$ board, which contains in turn
two cp. To prove all these facts, one can use coordinates. We
refer to Section \ref{s:patterns} for more details.  

The strategy is now clear. We shall study the structure of the
elementary boards, which are obtained in \cite{DeM} and look for
specific patterns into them. Then, we shall conclude by induction on the
dimension. We obtain: 

\begin{theorem}\label{t:main}
Let $2\leq n_1\leq n_2\leq... \leq n_k$, with $k\geq 3$. The
$n_1\times \ldots \times n_k$ chessboard 
 has \emph{no} closed knight tour if and only if
one of the following assumption holds:
\begin{enumerate}
\item For all $i$, $n_i$ is odd,
\item $n_{k-1}=2$,
\item $n_k=3$.
\end{enumerate} 
\end{theorem}
Note that the hypotheses are the same as the ones given in Theorem
\ref{t:DeM} when $k=3$. In the same paper they asked about higher
dimensional tours, This question was also asked by DeMaio \cite{De}
and Watkins \cite{Wat3}.  We mention that a conjecture for this
theorem was given in \cite{Kum}.    

The paper is organized as follows. In Section \ref{s:notation} we set
the notation and define the graph structure induced by a knight on an
$n$ dimensional board. Then, in Section \ref{s:patterns} we introduce
the notion of \emph{cross-patterns} and explain how to use them in
order to gain a dimension.
In Section \ref{s:forbidden}, we start the
proof of Theorem \ref{t:main} and provide all the negative
answers.  After that, in Section \ref{s:main}, we finish the proof of
Theorem \ref{t:main}. In Section \ref{s:conj} we show we can apply
this technique to  the problem of knight's tours with more general
moves. Finally, in the appendix,  
we give an alternative proof of Theorem \ref{t:DeM}.  

\noindent
{\bf Acknowledgments: } We would like to thank Imre Leader for useful
discussions and comments on the script. 

\section{Notation}\label{s:notation}
Let $\underline n:=(n_1, \ldots, n_k)$ be a multi-index where $n_i$
are with values in $\N\setminus\{0\}$. We denote by $|\underline
n|:=k$, the size of the multi-index $\underline n$. Set also
$\lint a,b\rint:= [a,b]\cap \Z$. The chessboard associated to
$\underline n$ is defined by 
\[\Bc_{{\underline n}}:= \lint 1, n_1\rint\times \ldots \times \lint
1, n_{|\underline n|}\rint.\] 
We turn to the definition of the
\emph{moves of the knight}. We set: 
\begin{align*}
\Cc_{|{\underline n}|}:=&\{(a_1, \ldots,  a_{|\underline n|})\in
\Z^{|{\underline n}|}, \quad  \mbox{ such that } 
\\ 
&\hspace{-0.5cm}|\{i, a_i=0\}|= |{\underline n}| -2,\quad  |\{i, a_i\in \{\pm 1\}\}|= 1, \quad  \mbox{ and } \quad  |\{i, a_i\in \{\pm 2\}\}|= 1\}.
\end{align*}

We endow $\Bc_{{\underline n}}$ with a graph structure, as follows.
We set $\Ec_{\underline   n}:\Bc_{{\underline n}} \times
\Bc_{{\underline n}} \to \{0,1\}$ be a symmetric function defined as follows: 
\begin{align}\label{e:graphrela}
\Ec_{\underline n} (a,b) := 1, \quad \mbox{ if } \quad a-b:= (a_1 - b_1,
\ldots, a_{|\underline n|}- b_{|\underline n|})\in
\Cc_{|{\underline n}|}
\end{align}
and $0$ otherwise. In other words, $a$ is linked to $b$ by a knight
move, if and only if $\Ec_{\underline n} (a,b) := 1$. The couple
$\Gc_{\underline n}:=(\Bc_{\underline n}, \Ec_{\underline n})$ is the
\emph{non-oriented graph} corresponding to all the possible paths of a
knight on the chessboard associated to $\underline n$. 

Set $\phi:\Bc_{\underline n}\to \{-1,1\}$ given by $\phi(a_1, \ldots, a_{|\underline n|}):=
(-1)^{a_1\times \ldots   \times a_{|\underline n|}}$, i.e., we assign the color black
or white to each square. Then note that given $a,b\in
\Bc_{\underline n}$ such that $\Ec_{\underline n}(a,b)=1$, one has $\phi(a)\times
\phi(b)=-1$. Therefore the graph $\Gc_{\underline n}:=(\Bc_{\underline
  n}, \Ec_{\underline n})$ is 
\emph{bipartite}.  
 
Let $\{a^i\}_{i=1, \ldots, |\Bc_{\underline n}|}$ be some elements of
$\Bc_{\underline n}$. We say that $\{a^i\}_{i}$ is a \emph{Hamiltonian cycle}
if the elements are two by two distinct and if $\Ec_{\underline n}
(a^i,a^{i+1}) := 1$ for all $i\in \{1, \ldots, |\Bc_{\underline n}|\}$
and if $\Ec_{\underline n} (a^1,a^{|\Bc_{\underline n}|}) := 1$. Note
that because of \eqref{e:graphrela}, $\Bc_{\underline n}$ has a closed knight tour
is just rephrasing that fact that $\Gc_{\underline n}$ has a Hamiltonian cycle. 

\begin{remark}\label{r:permu}
Given a multi-index $\underline{n}$ and $\psi$ a bijection from $\lint
1, |\underline{n}|\rint$ onto itself, we set $\underline{m}$ of size
$|\underline{n}|$ given by  $m_i:=n_{\psi(i)}$ for all $i\in \lint 1,
|n|\rint$. Since $\Cc_{|\underline{n}|}= \Cc_{|\underline{m}|}$ is
invariant under permutation,  $\Gc_{\underline n}$ has a Hamiltonian
cycle if and only $\Gc_{\underline m}$ has one. 
\end{remark}

\section{Looking for patterns}\label{s:patterns}
In Proposition \ref{p:gain}, see also Example \ref{ex:2d} for
  $2$-dimensional boards, we shall explain how to gain dimensions, for
  the question of closed tours,  with the help of patterns.

\begin{figure}[ht]
\centering
\includegraphics[scale=0.2]{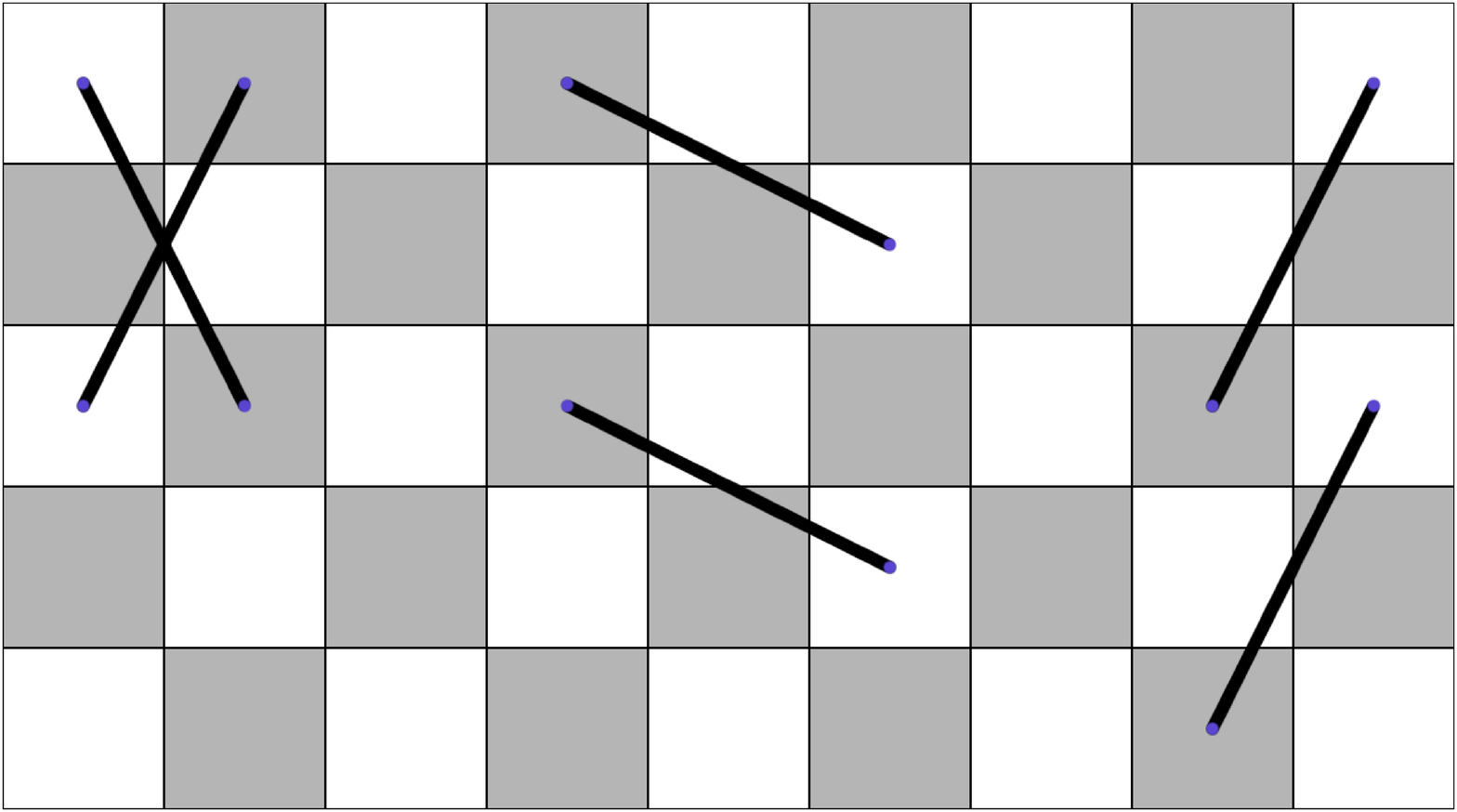}
\caption{Sites in dimension $2$.}\label{f:sites}
\end{figure}

Given an $n \times m$ board, we say that a pair of edges
$((a_1,b_1),(a_2,b_2))$ and   $((c_1,d_1),(c_2,d_2))$ in a tour is a
\emph{site} if  
\[(|a_1-c_1|,|b_1-d_1|)= (|a_2-c_2|,|b_2-d_2|)
\in \{ (0,2),(2,0)\}\]
or if
\[(|a_1-c_2|,|b_1-d_2|)= (|a_2-c_1|,|b_2-d_1|)\in \{
(0,2),(2,0)\}.\]
Roughly speaking, both endpoints of the two edges
are two squares away from each other. The three possible
configurations are given in Figure \ref{f:sites}. The first will be
\emph{cross-pattern} and the two others \emph{parallel-patterns}. 

\begin{figure}[ht]
\centering
\includegraphics[scale=0.2]{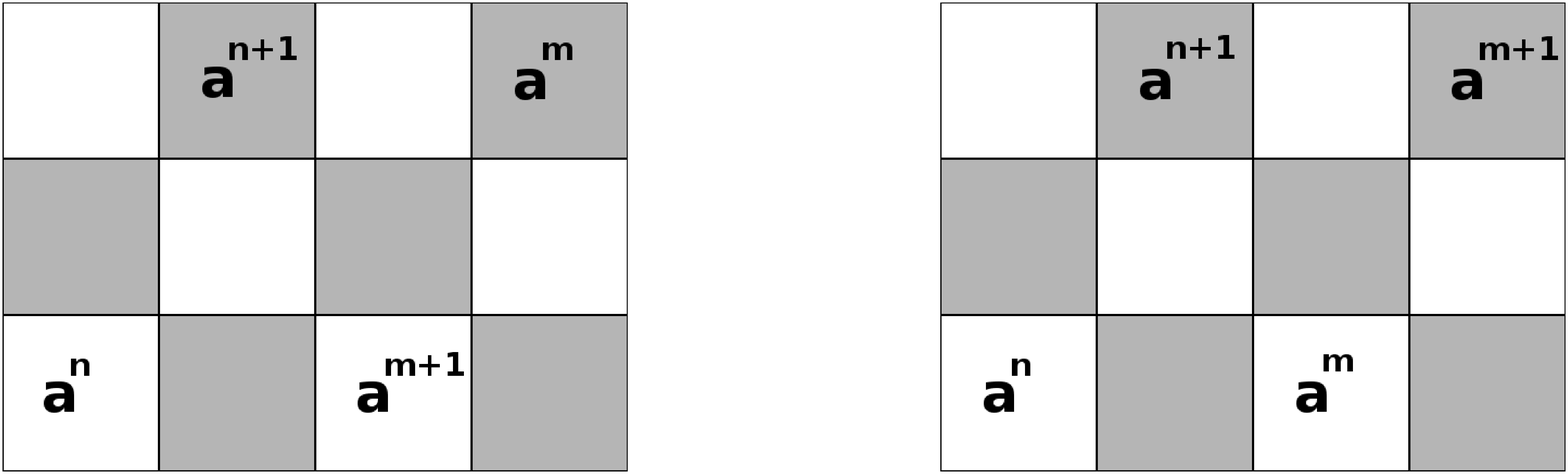}
\caption{A wopp and a nwopp.}\label{f:wopp}
\end{figure}

We generalize the definition for higher dimensional board. We
  denote by $(e_i)_{i\in \lint 1,   d\rint}$ the canonical basis of
  $\R^d$.
\begin{definition} 
Given a Hamiltonian cycle $(a^i)_{i\in I}$ for $\Gc_{\underline n}$, we say that it
contains the \emph{well-oriented parallel-pattern} (wopp) if there are
$c\in\Cc_{|{\underline n}|}$, $n, m \in I$, and $i\in \lint 1,
  d\rint$ so that:
\begin{align}\nonumber
a^{n+1}=a^n+c & \mbox { and } a^{m}= a^{m+1}+ c
\\
\label{e:wopp}
a^m- a^{n+1}=a^{m+1}-a^n &\in \{ \pm 2e_i\}.
\end{align}
We denote this wopp by $(a^n, a^{n+1}, a^m, a^{m+1})$. 

We say that it contains the \emph{non-well-oriented parallel-pattern}
(nwopp) if there are  $c\in\Cc_{|{\underline n}|}$, $n, m \in I$, and
$i\in \lint 1,
  d\rint $ so that: 
\begin{align}\nonumber
a^{n+1}=a^n+c & \mbox { and } a^{m+1}= a^{m}+ c
\\
\label{e:nwopp}
a^{m}-a^n = a^{m+1}-a^{n+1}&\in \{\pm 2e_i\}
\end{align}
We denote this nwocp by $(a^n, a^{n+1}, a^{m+1}, a^m)$. 
\end{definition}
Note that the equality in \eqref{e:wopp} (resp.\ in \eqref{e:nwopp}) is
automatically satisfied by the condition with $c$. We give an example
in a $2$-dimensional board in Figure \ref{f:wopp}.

We turn to cross-patterns. To each element $c:= (a_1,  
\ldots,  a_{\underline n})\in \Cc_{|{\underline n}|}$, we associate
two elements $e_{[c, 1]}$ and $e_{[c, 2]}$ of $\Z^{|{\underline n}|}$, 
by setting
\[e_{[c, 1]}(i):= \delta_{a_i\in \{\pm 1\}} \mbox{ and } e_{[c,
  2]}(i):= \delta_{a_i\in \{\pm 2\}},\] 
where $\delta_{X}:=0$ if $X$ is empty and $\delta_{X}:=1$ otherwise.
Note that: 
\[c = \langle c, e_{[c, 1]}\rangle e_{[c, 1]} + \langle c,
e_{[c, 2]}\rangle e_{[c, 2]}\in \{\pm 1 e_{[c, 1]}+ \pm 2 e_{[c,  2]}\},\]
where $\langle \cdot, \cdot\rangle$ denotes the Euclidean scalar
product in $\R^{|\underline{n}|}$. 
We now flip $c$ with respect to $e_{[c, 1]}$. More precisely, we set
\[\tilde c:= -\langle c, e_{[c, 1]}\rangle e_{[c, 1]} + \langle c,
e_{[c, 2]}\rangle e_{[c, 2]}.\] 

\begin{definition} 
Given a Hamiltonian cycle $(a^i)_{i\in I}$ for $\Gc_{\underline n}$, we say that it
contains the \emph{well-oriented cross-pattern} (wocp) if there are
$c\in\Cc_{|{\underline n}|}$ 
and  $n, m \in I$, so that:
\begin{align}\label{e:wocp}
a^{n+1}=a^n+c & \mbox { and } a^{m+1}= a^m+\tilde c
\\
\nonumber
a^m-a^n &= \langle c, e_{[c, 1]}\rangle  e_{[c, 1]}.
\end{align}
We denote this wocp by $(a^n, a^{n+1}, a^m, a^{m+1})$. 

We say that it contains the \emph{non-well-oriented cross-pattern}
(nwocp) if there are  $c\in\Cc_{|{\underline n}|}$ 
and  $n, m \in I$, so that:
\begin{align}\label{e:nwocp}
a^{n+1}=a^n+c & \mbox { and } a^{m}= a^{m+1}+\tilde c
\\
\nonumber
a^{m+1}-a^n &=  \langle c, e_{[c, 1]}\rangle e_{[c, 1]}.
\end{align}
We denote this nwocp by $(a^n, a^{n+1}, a^{m+1}, a^m)$. 
\end{definition}

\begin{figure}[ht]
\centering
\includegraphics[scale=0.2]{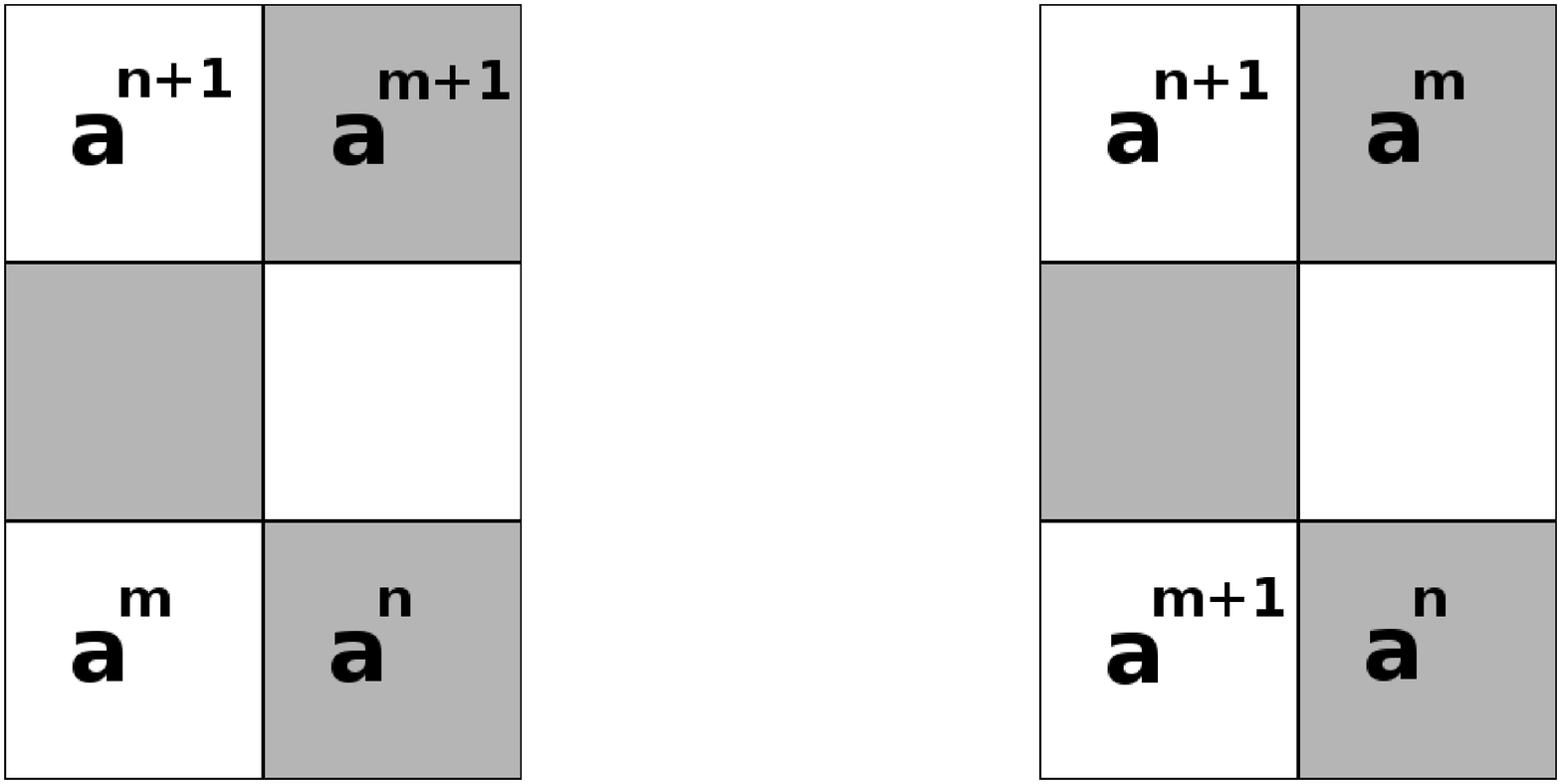}
\caption{A wocp and a nwocp.}\label{f:wocp}
\end{figure}

\begin{definition}
 Given a Hamiltonian cycle $(a^i)_{i\in I}$, we say that $(a^n,
 a^{n+1}, a^{m}, a^{m+1})$ is a well-oriented site (or pattern) if it is
 a wocp or a wopp. we say that $(a^n,
 a^{n+1}, a^{m+1}, a^m)$ is a non-well-oriented site (or pattern) if it is
 a nwocp or a nwopp.
\end{definition} 

The Properties \eqref{e:wopp} and \eqref{e:nwopp} are important and we
present the equivalent for cross-patterns. It would be
extensively used  in Example \ref{ex:2d} and in the Proposition \ref{p:gain}.

\begin{remark}\label{r:key} In the case of the wocp \eqref{e:wocp}, we have 
\[a^{n+1}-a^m=a^{m+1}- a^n=\langle c,
e_{[c, 2]}\rangle e_{[c, 2]}.\]
In particular, if $(a^n, a^{n+1}, a^{m+1}, a^m)$ is a well-oriented
site, then there is an $i$ such that \[a^{n+1}-a^m= \pm (a^{m+1}-a^n)
\in \{-2 e_i, 2e_i\}.\] 
  
In the case of the nwocp \eqref{e:nwocp}, we get
\[a^{n+1}-a^{m+1}=a^{m}- a^n=\langle c,
e_{[c, 2]}\rangle e_{[c, 2]}.\]
This also yields that, if $(a^n, a^{n+1}, a^{m+1}, a^m)$ is a
non-well-oriented site, then there is an $i$ such
that \[a^{n+1}-a^{m+1}=\pm(a^{m}-a^n) \in \{-2 e_i, 2e_i\}.\]
\end{remark}

We now explain how to connect two Hamiltonian cycles and gain one
dimension with the help of sites. We mention that, in the
$2$-dimensional case, the existence of sites will be discussed in
Proposition \ref{p:2dbi}.  
\begin{example}\label{ex:2d}
Take $\Gc_{\underline
  {n}}$, with $|\underline
  {n}|=2$, such that it contains a Hamiltonian cycle
$(a^i)_{i\in I}$, see Theorem \ref{t:S}. Take now
$\underline{m}:=(\underline{n}, 2)$. Note that $\Gc_{\underline{m}}$ contains two
cycles:
\[ (a^i, 1)_{i\in I} \mbox{ and } (a^i, 2)_{i\in I}.
\] 
Suppose there is a well-oriented site $(a^n, a^{n+1}, a^m,
a^{m+1})$. Using Remark \eqref{r:key}, we obtain immediately that
$\Ec_{\underline{m}}((a^{n}, 1),(a^{m+1}, 2))=1$ 
and $\Ec_{\underline{m}}((a^{m}, 2),(a^{n+1}, 1))=1$. 
We can construct a Hamiltonian cycle as follows:
\begin{align*}
(a^{n}, 1)\to (a^{m+1}, 2) \to (a^{m+2}, 2) \to \ldots \to (a^{m}, 2)
\to (a^{n+1}, 1) \to \ldots \to (a^{n}, 1).
\end{align*}
Suppose now there is a non-well-oriented site $(a^n, a^{n+1}, a^{m+1},
a^{m})$. Similarly, We construct a Hamiltonian cycle as follows: 
\begin{align*}
(a^{n}, 1)\to (a^{m}, 2) \to (a^{m-1}, 2) \to \ldots \to (a^{m+1}, 2)
\to (a^{n+1}, 1) \to \ldots \to (a^{n}, 1).
\end{align*}
Note that we went backward on the second copy. 
\end{example}

It is obvious that one can gain as many dimension as one has disjoint
cross-patterns. However, the situation is much better: having solely
two of them are enough to gain as many dimension as we want.

\begin{definition}
 A tour would be called \emph{bi-sited} if it contains two sites with
disjoint support,  i.e., such
that the endpoints of the four pair of edges are two by two disjoint.
\end{definition} 

\begin{figure}[ht]
\centering
\includegraphics[scale=0.2]{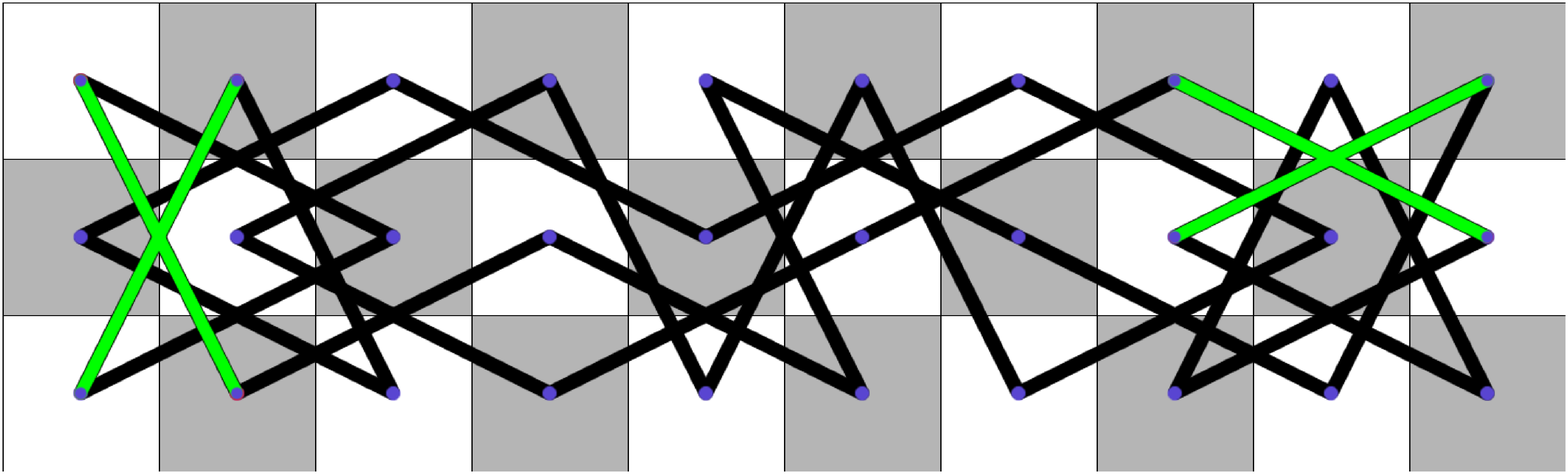}
\caption{A bi-sited closed tour}\label{f:bi-sited} 
\end{figure}

We rely on the next Proposition.  

\begin{proposition}\label{p:gain}
Take $\Gc_{\underline
  {n}}$. Suppose that it contains a bi-sited Hamiltonian cycle. Take
now $\underline{m}:=(\underline{n}, k)$, with $k\geq 2$. Then 
$\Gc_{\underline{m}}$ contains a bi-sited Hamiltonian cycle.
\end{proposition}

\proof As the demonstration is similar, we shall present only the case
of two well-oriented sites. Set $(a^i)_{i\in I}$ the Hamiltonian
cycle. Note that $(a^i, j)_{i\in I}$, for $j\in \lint 1, k\rint$, are
$k$ disjoint cycle in $\Gc_{\underline{m}}$. Take two well-oriented
sites  $(a^{n_1}, a^{n_1+1}, a^{m_1}, a^{m_1+1})$ and
$(a^{n_2}, a^{n_2+1}, 
a^{m_2}, a^{m_2+1})$. Let $k$ be \emph{even}. Using Remark
\ref{r:key}, we obtain 
\begin{align}
\nonumber
(a^{n_1}, 1)&\to (a^{m_1+1}, 2) \to (a^{m_1+2}, 2) \to \ldots \to
(a^{n_2}, 2) \to
\\
\nonumber
& \to (a^{m_2+1}, 3) \to (a^{m_2+2}, 3)\to \ldots \to (a^{n_1}, 3)\to 
\\
\nonumber
& \to (a^{m_1+1}, 4) \to (a^{m_1+2}, 4) \to \ldots \to (a^{n_2}, 4)
\to \ldots
\\
\label{e:gain1}
& \to (a^{m_2+1}, k-1) \to (a^{m_2+2}, k-1) \to \ldots \to (a^{n_1}, k-1) \to
\\
\label{e:gain2}
& \to (a^{m_1+1}, k) \to (a^{m_1+2}, k) \to \ldots \to (a^{m_1}, k)
\to
\\
\label{e:gain3}
& \to (a^{n_1+1}, k-1) \to (a^{n_1+2}, k-1) \to \ldots \to (a^{m_2},
k-1) \to \ldots
\\
\nonumber
& \to (a^{n_2+1}, 2) \to (a^{n_2+2}, 2) \to \ldots \to (a^{m_1}, 2)
\\
\nonumber
&\to  (a^{n_1+1}, 1)\to (a^{n_1+2}, 1)\to \ldots \to (a^{n_1}, 1).
\end{align}
Finally note that $((a^{n_2}, j), (a^{n_2+1},j), (a^{m_2},j),
(a^{m_2+1}, j))$ with $j\in\{1, k\}$ are well-oriented sites
for this Hamiltonian cycle in $\Gc_{\underline m}$. When $k$ is odd,
we replace the lines \eqref{e:gain1}, \eqref{e:gain2}, and
\eqref{e:gain3} by
\begin{align*}
\quad\quad\quad \mbox{ }
&\to (a^{m_1+1}, k-1) \to (a^{m_1+2}, k-1) \to \ldots
\to (a^{n_2}, k-1) \to 
\\
& \to (a^{m_2+1}, k) \to (a^{m_2+2}, k) \to \ldots \to (a^{m_2}, k)
\to
\\
& \to (a^{n_2+1}, k-1) \to (a^{n_2+2}, k-1) \to \ldots \to (a^{m_1},
k-1) \to \ldots
\end{align*}
 In this case, we get
$((a^{n_2}, 1), (a^{n_2+1},1), (a^{m_2},1), 
(a^{m_2+1}, 1))$ and  
$((a^{n_1}, k), (a^{n_1+1},k), (a^{m_1},k), 
(a^{m_1+1}, k))$  are well-oriented sites associated to the
Hamiltonian cycle in $\Gc_{\underline m}$. \qed

\section{Into the proof}
In this section we give the proof of Theorem \ref{t:main} and consider
boards associated to a multi-index $\underline n:=(n_1, \ldots,
n_{|\underline n|})$, where $n_i\geq 2$ and $|\underline{n}|\geq 3$.  
\subsection{Forbidden boards}\label{s:forbidden} We start with the
sufficient part of Theorem \ref{t:main}. It is a straightforward
generalization of the dimension $3$. There are three main cases: 

a) $n_i\in (2\N+1)$, for all $i$.  A bipartite graph cannot have a
Hamiltonian cycle if its cardinality is odd.  

b) $\underline n:=(2, 2, \ldots, 2, k)$, with $k\geq 2$. Set $a^0:=(2, 2,
\ldots, 2)\in \Bc_n$. Let $(a^i)_{i\geq 0}$, be the elements of a
Hamiltonian cycle. Note that $(a^0-a^i)\in (\{1,2\}, \ldots,
\{1,2\}, 2\Z)$, for all $i\geq 0$. Then, $(1, 1, \ldots, 1)$ is never
reached. Contradiction.  

c) $\underline n:=(n_1, \ldots, n_{|\underline n|})$, with 
$2\leq n_1\leq \ldots \leq n_{|\underline{n}|}\leq 3$. Set $a:=(2,
\ldots, 2)\in \Bc_n$ and note that there is no $b\in
\Cc_{|\underline{n}|}$, such that $a+b\in
\Bc_{\underline{n}}$. Therefore, the graph
$\Gc_{\underline{n}}:=(\Bc_{\underline{n}}, 
\Ec_{\underline{n}})$ has no closed tour. 

\begin{remark}\label{r:disco}
Note that the two last cases correspond to disconnected boards.
\end{remark} 

\subsection{Bi-sited boards in dimension $2$ and  $3$}
We start by the existence in dimension $2$. 
\begin{proposition}\label{p:2dbi}
Every closed tour on an $n\times m$ board is bi-sited.
\end{proposition} 
\proof First, there is no closed tour on a $2\times k$ and on $4\times
k$ boards.  

\noindent {\bf i) Case ${\bf 3\times k}$, for ${\bf k\geq 10}$:} Note that
$k$ is even. First one has $(1,1)$ which is linked to $(3,2)$ and to
$(2,3)$. The possible neighbors of $(1,3)$ are  $(3,2)$,$(2,1)$,
$(3,4)$, $(2,5)$. As we have a cycle, exactly two of them are linked
to $(1,3)$ and three of them are part of a pattern (the second one
forms a cross-pattern and the two last one parallel patterns). Then,
there is at least one site. Secondly we repeat the argument for the upper
right corner and get a new site. Since $k\geq 10$, the two sites are
with disjoint supports. 

\noindent {\bf ii) Case ${\bf k\times l}$, for ${\bf k\geq 5,}$  and
  ${\bf l\geq 6}$:} We flip the board and consider $l\times k$. We
repeat twice the first part of the point i) for the sub-board of size
$3\times k$ that contains the upper left corner and  for the one that
contains the lower left corner. We get two disjoint sites. \qed

This proposition will be very useful in Appendix
\ref{s:DeMalt}. Indeed, combining it with Proposition \ref{p:gain}, we infer:

\begin{corollary}\label{c:5}
If an $n \times m$ closed tour exists then so does an $n \times m \times
p_1 ....\times p_r$ closed tour for any $p_1, \ldots, p_r \in
\N\setminus \{0\}$.
\end{corollary} 

We were not able to obtain an analogue to Proposition \ref{p:2dbi} for $3$ dimensional boards.
 However, for our purpose it is enough to prove the existence of specific
bi-sited tours. In this section
we will rely on the construction of \cite{DeM} and in Appendix
\ref{s:DeMalt}, we will give a self-contained proof.

\begin{theorem}\label{t:DeM2} Let $2\leq m\leq n\leq
  p$.  The $m\times   n\times p$  chessboard has a bi-sited closed
  tour if and only  if one of the following assumption holds:
\begin{enumerate}
\item $m$, $n$, or $p$ is even,
\item $n\geq 3$,
\item $p\geq 4$.
\end{enumerate} 
\end{theorem}
\proof Without going into details, we give a rough idea about the approach
 of \cite{DeM}. Take $n\times m\times p$, that is not satisfying the
 hypotheses a), b) and c) of Theorem \ref{t:DeM}. Note that at least 
one of them is even (say $n$). 
Subsequently, one basically works modulus $4$: An $n\times m\times p$
block can be written as a union of the following ones: $2\times 4\times 4$,
$2\times 4\times 5$, $2\times 4\times 6$, $2\times 4\times 3$,
$2\times 5\times 5$, $ 2\times 5\times 6$, $2\times 5\times 3$,
$2\times 6\times 6$, $2\times 6\times 3$,  $2\times 7\times 3$,
$4\times 3\times 3$, $6\times 3\times 3$. The list is rather long
since one has no Hamiltonian cycle for $2\times 2\times 3$ and
$2\times 3\times 3$ boards. Then, the authors
construct some Hamiltonian cycle for each elementary blocks. 
After that, given two compatible blocks, by deleting one edge from each
block and by creating two edges that are ``gluing'' these blocks
together, they  construct a Hamiltonian cycle for the union, starting
with the two disjoint ones. 

First, we study the elementary blocks that are exhibited
in \cite{DeM}. For each of them, we give the list of well-oriented
cross-patterns (wocp), non-well-oriented cross-patterns (nwocp) and
also the edges that are used to combine two different blocks, in order
to create a Hamiltonian cycle for the union. We will not exhibit 
parallel-patterns. For the latter, we indicate the other figure that is
glued with and denote it between square brackets. We also strike out
all the wocp and nwocp which are incompatible with the gluing
operation of \cite{DeM}. For instance, in Figure $6$, we struck
$(2,3,21,20)$ out, because the edge $(20,21)$ is already used when one
glues Figure $35$.

\begin{align*}
\begin{array}{c}
\begin{array}{|c|c|c|}
\hline
\mbox{Wocp}&\mbox{Nwocp}&\mbox{Used edges [with Figure]}
\\
\hline
(4,5,11,12)&(1,2,22,21)& (3,4)[6,11,14,17],(8,9)[6]
\\
(7,8,28,29)&\hcancel{(2,3,21,20)}& (10,11)[6],(14,15)[6],
\\
&(13,14,18,17)& (20,21) [35], (25,26)[6],
\\
&&(27,28)[6]
\\
\hline
\end{array}
\\
\\
\mbox{Case of the Figure } 6: \mbox{ size } 2\times 4\times 4
\end{array}
\end{align*}

\begin{align*}
\begin{array}{c}
\begin{array}{|c|c|c|}
\hline
\mbox{Wocp}&\mbox{Nwocp}&\mbox{Used edges [with Figure]}
\\
\hline
(1,2,30,31)&(14,15, 35,34)& (4,5) [20],
(8,9) [20], 
\\
\hcancel{(4,5,31,32)}&\hcancel{(16,17,23,22)}& (10,11)[11], 
(11,12) [23],
\\
\hcancel{(12,13,39,40)}&&(16,17)[11], 
(17,18) [26], 
\\
(19,20,24,25)&& (31,32)[11], (37,38)[11], (39,40)[6]
\\
\hline
\end{array}\\
\\
\mbox{Case of the Figure } 11: \mbox{ size } 2\times 4\times 5
\end{array}
\end{align*}

\begin{align*}
\begin{array}{c}
\begin{array}{|c|c|c|}
\hline
\mbox{Wocp}&\mbox{Nwocp}&\mbox{Used edges [with Figure]}
\\
\hline
(1,2,3,4)&(14,15, 47, 46)& 
(4,5) [29], 
(9,10) [32], 
\\
(7,8,26,27)&&(28,29)[14],(29,30)[14],
\\
(8,9,25,26)&&(30,31)[14],(36,37)[6], 
\\
(10,11, 31,32)&&(39,40)[14],(43,44) [29],
\\
(11,12, 44,45)&&(47,48) [23]
\\
(13,14,32,33)&&
\\
(15,16,34, 35)&&
\\
\hcancel{(18,19, 39,40)}&&
\\
(48,1,37,38)&&
\\
\hline
\end{array}
\\
\\
\mbox{Case of the Figure } 14: \mbox{ size } 2\times 4\times 6
\end{array}
\end{align*}
\begin{align*}
\begin{array}{c}
\begin{array}{|c|c|c|}
\hline
\mbox{Wocp}&\mbox{Nwocp}&\mbox{Used edges [with Figure]}
\\
\hline
(8,9,21,22)&(2,3,7,6)& (4,5)[17],(5,6)[17], 
\\
(9,10,16,17)&&(7,8)[17],(11,12)[6],   
\\
(10,11, 15,16)&& (20,21)[17, 26, 32, 35]
\\
\hline
\end{array}
\\
\\
\mbox{Case of the Figure } 17: \mbox{ size } 2\times 4\times 3
\end{array}
\end{align*}
\begin{align*}
\begin{array}{c}
\begin{array}{|c|c|c|}
\hline
\mbox{Wocp}&\mbox{Nwocp}&\mbox{Used edges [with Figure]}
\\
\hline
(3,4,38,39)&(9,10, 24, 23)& (20,21)[11], (43,44)[20],
\\
(5,6,28,29)&(18,19, 33,32)&(47,48)[11],  (49,50)[20]
\\
(8,9,25,26)&\hcancel{(19,20,50,49)}&
\\
(12,13,45,46)&\hcancel{(20,21,43,42)}&
\\
\hcancel{(14,15,47,48)}&&
\\
(17,18,30,31)&&
\\
\hline
\end{array}
\\
\\
\mbox{Case of the Figure } 20: \mbox{ size } 2\times 5\times 5
\end{array}
\end{align*}

One has to be careful with the block $2\times 5 \times 6$
which is given in Figure $23$ of \cite{DeM}. Indeed, one sees that $1$
cannot be reached from $60$, therefore this is not a Hamiltonian cycle
but just a path. To fix this, we propose to take:

\begin{figure}[ht]
\centering
\includegraphics[scale=0.2]{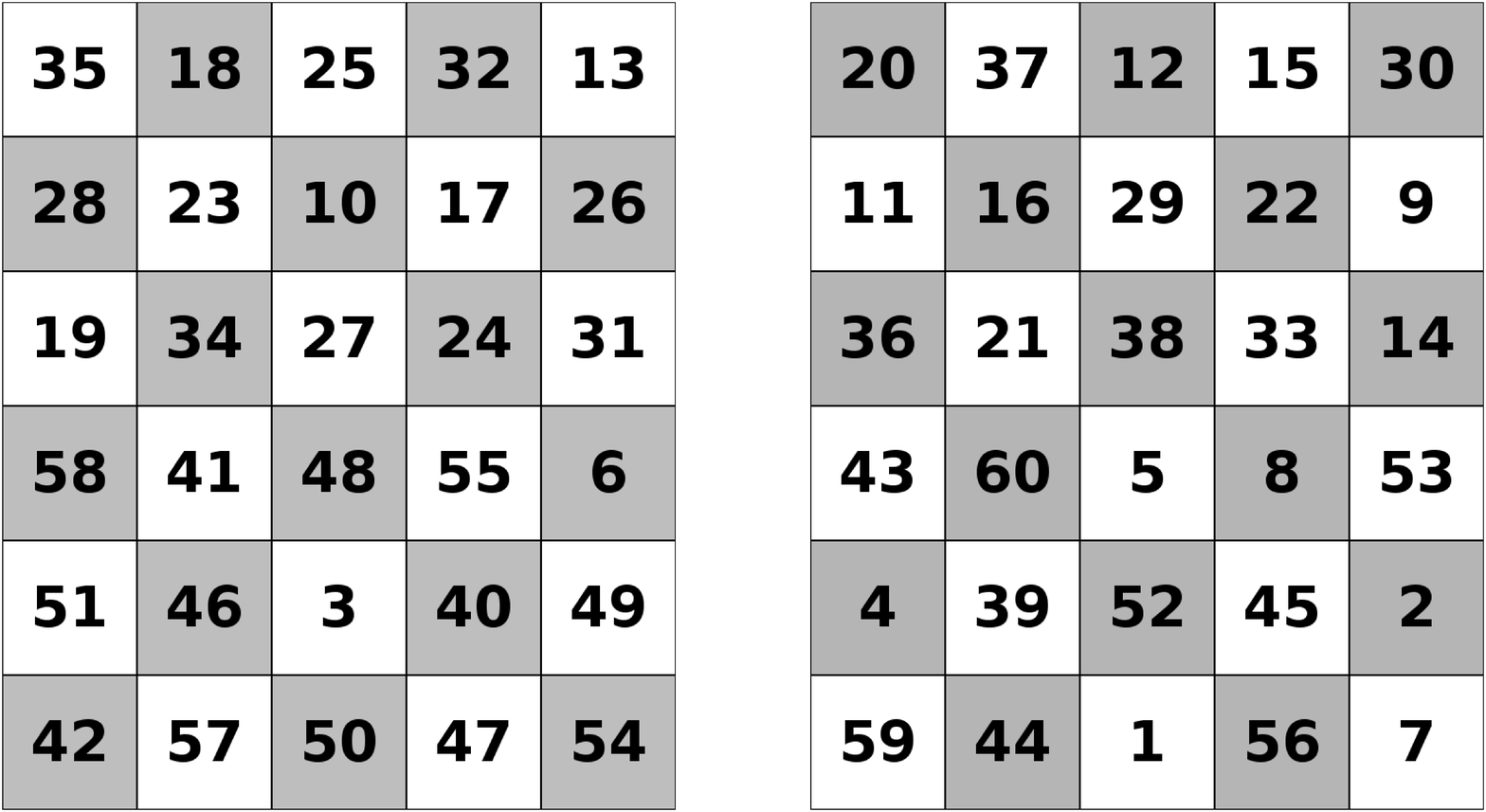}
\caption{A $5\times 6\times 2$ bi-sited open tour replacing Figure 23
  of \cite{DeM}.}
\end{figure}

We  replace the table following the board, i.e., the
table on the top of page $9$. By 

\begin{center}
\begin{tabular}{ |l | l | l|}
\hline
&Delete edges& Create edges
\\
\hline
Vertical& $11-12$ top board, $12-13$ bottom board&$11-12$, $12-13$
\\
\hline
Horizontal&  $32-33$ left board, $47-48$ right board& $32-47$, $33-48$
\\
\hline
Front& $1-2$ front board, $44-45$ back board& $1-44$, $2-45$
\\
\hline
\end{tabular}
\end{center}

This yields:
\begin{align*}
\begin{array}{c}
\begin{array}{|c|c|c|}
\hline
\mbox{Wocp}&\mbox{Nwocp}&\mbox{Used edges [with Figure]}
\\
\hline
(6,7,53,54)&(3,4,52,51)&(1,2)[23], (12,13)[11],
\\
(13,14,30,31)&(10,11,29,28)&(32,33)[14], (44,45)[23]
\\
&(16,17,23,22)&
\\
&(18,19,35,34)&
\\
&(19,20,36,35)&
\\
&(20,21,37,36)&
\\
&(39,40, 46,45)&
\\
&(41,42,58,57)&
\\
&(42,43,59,58)&
\\
&(43,44,60,59)&
\\
\hline
\end{array}
\\
\\
\mbox{Case of the new Figure } 23: \mbox{ size } 2\times 5\times 6
\end{array}
\end{align*}
\begin{align*}
\begin{array}{c}
\begin{array}{|c|c|c|}
\hline
\mbox{Wocp}&\mbox{Nwocp}&\mbox{Used edges [with Figure]}
\\
\hline
\hcancel{(6,7,23,24)}&(3,4,22,21)& (14,15)[11], (23, 24)[17], 
\\
&(9,10,16,15)&(27,28)[26], (29,30)[26]
\\
&\hcancel{(11,12,28,27)}&
\\
&(12,13,29,28)&
\\
&\hcancel{(13,14,30,29)}&
\\
\hline
\end{array}
\\
\\
\mbox{Case of the Figure } 26: \mbox{ size } 2\times 5\times 3
\end{array}
\end{align*}
\begin{align*}
\begin{array}{c}
\begin{array}{|c|c|c|}
\hline
\mbox{Wocp}&\mbox{Nwocp}&\mbox{Used edges [with Figure]}
\\
\hline
\hcancel{(15,16,62,63)}&(2,3,69,68)& (15, 16)[14],  (22,23)[29], 
\\
(16,17,57,58)&(4,5,47,46)& (29, 30)[14], (32,33)[29]
\\
(47,48,72,1)&(6,7,37,36)&
\\
&(7,8,36,35)&
\\
&(10,11, 51,50)&
\\
&(12,13, 39,38)&
\\
&(20,21,53,52)&
\\
&(23,24, 70,69)&
\\
&(27,28, 64,63)&
\\
&(30,31, 43,42)&
\\
&(31,32, 42,41)&
\\
\hline
\end{array}
\\
\\
\mbox{Case of the Figure } 29: \mbox{ size } 2\times 6\times 6
\end{array}
\end{align*}
\begin{align*}
\begin{array}{c}
\begin{array}{|c|c|c|}
\hline
\mbox{Wocp}&\mbox{Nwocp}&\mbox{Used edges [with Figure]}
\\
\hline
(6,7,15,16)&(3,4,14,13)& (11,12)[14],  (18,19)[32], 
\\
(8,9,17,18)&(4,5,15,14)&(22,23)[32], (26,27)[17]
\\
\hcancel{(18,19, 31,32)}&(25,26,34,33)&
\\
(19,20, 28,29)&\hcancel{(26,27,35,34)}&
\\
\hline
\end{array}
\\
\\
\mbox{Case of the Figure } 32: \mbox{ size } 2\times 6\times 3
\end{array}
\end{align*}
\begin{align*}
\begin{array}{c}
\begin{array}{|c|c|c|}
\hline
\mbox{Wocp}&\mbox{Nwocp}&\mbox{Used edges [with Figure]}
\\
\hline
(2,3,25,26)& (3,4,24,23)& (7,8)[6],  (12,13)[35],
\\
(6,7,19,20)&(5,6,22,21)& (13,14)[17],(14,15)[35]
\\
\hcancel{(7,8,16,17)} &(32,33,39,38)&
\\
(8,9,15,16)&(33,34,40,39)&
\\
\hcancel{(13,14,37,38)}&&
\\
\hline
\end{array}
\\
\\
\mbox{Case of the Figure } 35: \mbox{ size } 2\times 7\times 3
\end{array}
\end{align*}
\begin{align*}
\begin{array}{c}
\begin{array}{|c|c|c|}
\hline
\mbox{Wocp}&\mbox{Nwocp}&\mbox{Used edges [with Figure]}
\\
\hline
&(3,4,36,35)& (9,10)[39], (28,29)[39]
\\
&\hcancel{(9,10,20,19)}&
\\
&(11,12, 16,15)&
\\
&\hcancel{(28,29, 1, 36)}&
\\
\hline
\end{array}
\\
\\
\mbox{Case of the Figure } 39: \mbox{ size } 4\times 3\times 3
\end{array}
\end{align*}
\begin{align*}
\begin{array}{c}
\begin{array}{|c|c|c|}
\hline
\mbox{Wocp}&\mbox{Nwocp}&\mbox{Used edges [with Figure]}
\\
\hline
(2,3,5,6)&(4,5,9,8)& (21,22)[40], (28,29)[40]
\\
(10,11, 15,16)&\hcancel{(18,19,29,28)}&
\\
(11,12,14,15)&(32,33,43,42)&
\\
(13,14,16,17)&&
\\
(19,20,22,23)&&
\\
(20,21,25,26)&&
\\
\hcancel{(21,22,24,25)}&&
\\
(36,37,49,50)&&
\\
\hline
\end{array}
\\
\\
\mbox{Case of the Figure } 40: \mbox{ size } 6\times 3\times 3
\end{array}
\end{align*}
In every case there is at least a pair cross-patterns which
have disjoint supports. Now thanks to the proof of Theorem
\ref{t:DeM}, we get there is a Hamiltonian cycle in $\Gc_{\underline
  n}$ which is obtained by gluing together Hamiltonian cycles coming from
the elementary blocks that we have just discussed. Therefore, the
Hamiltonian cycle is bi-sited.\qed

\subsection{Main result}\label{s:main} We finish the proof
of Theorem \ref{t:main}.

\proof[Proof of Theorem \ref{t:main}] The sufficient part was proved
in Section \ref{s:forbidden}.  We deal
with the necessary part by an induction on: 
\\
\\
$\Pc_k:=$``Given a multi-index $\underline{n}=(n_1, n_2, \ldots, n_{k})$, where 
$2\leq n_1\leq n_2\leq \ldots \leq n_k$ and such that:
\begin{enumerate}
\item There is $i_0\in\lint 1, k\rint$ such that $n_{i_0}$ is even,
\item $n_{k-1}\geq 3$,
\item $n_k\geq 4$,
\end{enumerate}
then, $\Gc_{\underline{n}}$ contains a bi-sited Hamiltonian
cycle."
\\
\\
for all $k\geq 3$.
\\

\noindent{\bf Basis:} Theorem \ref{t:DeM2} yields $\Pc_3$.  
\\

\noindent{\bf Inductive step:} Suppose that $\Pc_k$ holds true for some $k\geq
3$. Take $\underline n=(n_1, n_2, \ldots, n_{\bf k+1})$, where
$2\leq n_1\leq n_2\leq \ldots \leq n_{k+1}$, such that a), b) and c)
hold true.

If all $n_i$ but one are odd, then, up to a permutation, we can always
suppose that $n_1$ is odd. We set 
\begin{align}\label{e:n'}
\underline{n'}= (n_2, \ldots,n_{k+1}).  
\end{align}
Otherwise, without operating a permutation, we set \eqref{e:n'}.

Note that $\Gc_{\underline n'}$ satisfies the hypothesis a), b) and c)
of $\Pc_k$ and therefore, by $\Pc_k$, we get it contains a Hamiltonian
cycle with at least well-oriented or two disjoint  non-well-oriented
cross-patterns. We apply Proposition \ref{p:gain} and get
$\Pc_{k+1}$ is true. 

Therefore, we have proved by induction that $\Pc_k$ is true for all
$k\geq 3$. In particular, this proves the theorem. \qed

Using Remark \ref{r:disco}, we derive:

\begin{corollary}\label{c:7}
For $k\geq 3$ and $n_i\geq 2$. Set $\underline{n}:=(n_1, \ldots,
n_k)$. Suppose that some $n_i$ is even. Then, a closed tour exists on
$\Bc_{\underline{n}}$ if and only if $\Gc_{\underline{n}}$ is
connected.   
\end{corollary}

\section{Generalised knight's tours on a chessboard}\label{s:conj} 
The knight's tour is a specific case of many general questions. A
natural one to ask would be, what about more general moves? For
example instead of the knight being able to move $(\pm 1, \pm 2)$ or
$(\pm 2, \pm 1)$ what if the knight could move $(\pm \alpha, \pm \beta)$ or
$(\pm \beta, \pm \alpha)$?

Given a chessboard $\Bc_{{\underline n}}$ and $\alpha,\beta \in
\N\setminus\{0\}$ we define as before 
\begin{align*}
\Cc^{\alpha,\beta}_{|{\underline n}|}:=&\{(a_1, \ldots,  a_{|\underline n|})\in
\Z^{|{\underline n}|}, \quad  \mbox{ such that } 
\\ 
&\hspace{-0.5cm}|\{i, a_i=0\}|= |{\underline n}| -2,\quad  |\{i, a_i\in \{\pm \alpha\}\}|= 1, \quad  \mbox{ and } \quad  |\{i, a_i\in \{\pm \beta\}\}|= 1\}.
\end{align*}
and we endow $\Bc_{{\underline n}}$ with a graph structure, as follows.
We set $\Ec^{\alpha,\beta}_{\underline   n}:\Bc_{{\underline n}} \times
\Bc_{{\underline n}} \to \{0,1\}$ to be the symmetric function defined as follows: 
\begin{align*}
\Ec^{\alpha,\beta}_{\underline n} (a,b) := 1, \quad \mbox{ if } \quad a-b:= (a_1 - b_1,
\ldots, a_{|\underline n|}- b_{|\underline n|})\in
\Cc^{\alpha,\beta}_{|{\underline n}|}
\end{align*}
and $0$ otherwise. The couple
$\Gc^{\alpha,\beta}_{\underline n}:=(\Bc_{\underline n},
\Ec^{\alpha,\beta}_{\underline n})$ is the  
graph corresponding to all the possible paths of a generalised knight
on the chessboard associated  
to $\underline n$. Note that it is bipartite. 
We define an \emph{$(\alpha,\beta)$-tour} on
$\Bc_{\underline n}$ to be a Hamiltonian cycle on
$\Gc^{\alpha,\beta}_{\underline n}$. 

\begin{figure}[h!]\label{f:10x10}
\centering
\includegraphics[scale=0.2]{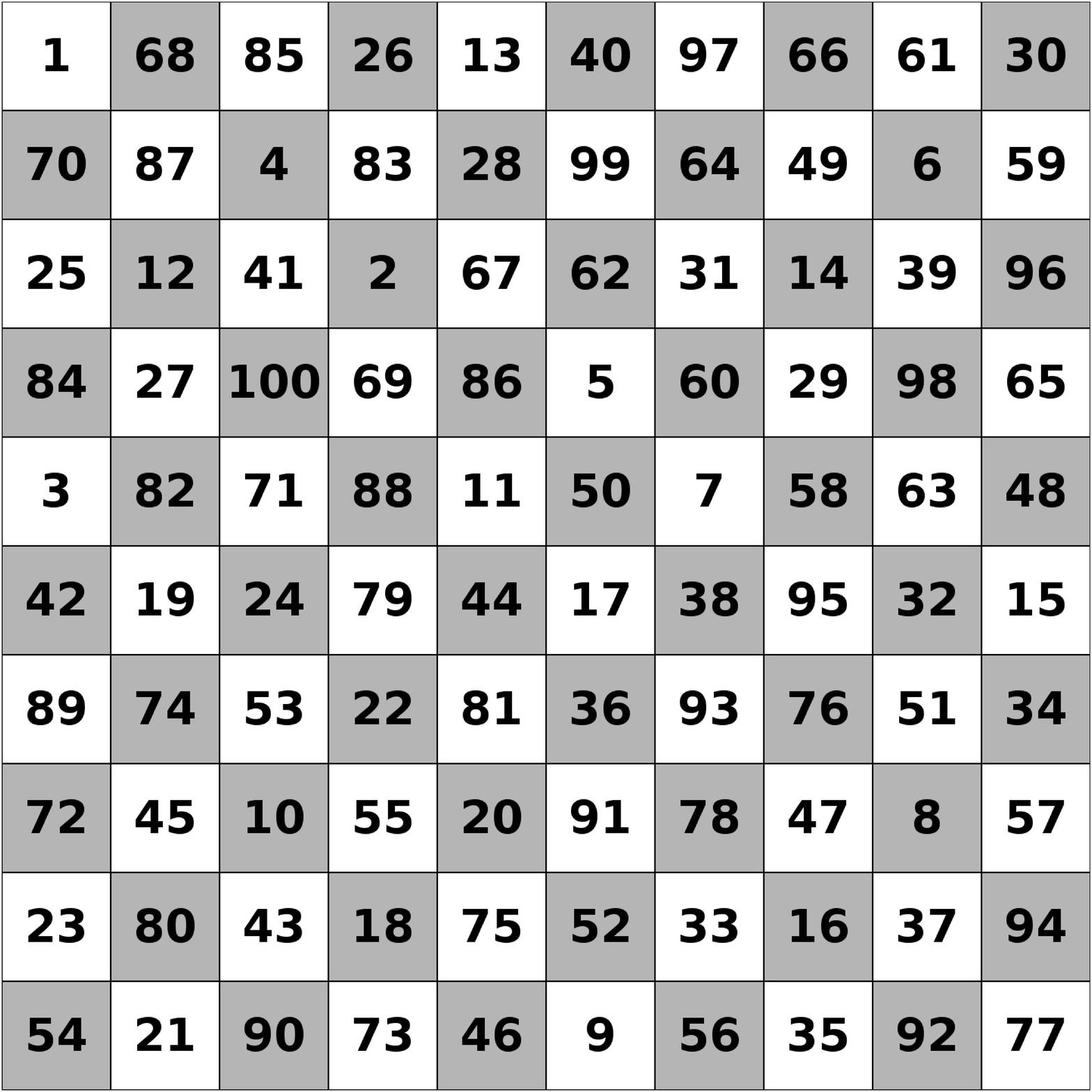}
\caption{A $10 \time 10$ $(2,3)$-tour constructed by A. H. Frost
\cite{Fro}}\label{f:10x10} 
\end{figure}

Knuth showed in \cite{Knu} that in $2$
dimensions if $\alpha \geq \beta$ then $\Gc^{\alpha,\beta}_{(n_1,n_2)}$
is connected if and only if $\gcd(\alpha+\beta,\alpha-\beta)=1$, $n_1 
\geq 2\alpha$, and $n_2 \geq \alpha+\beta$.

\begin{remark}\label{r:coprime}
If an $(\alpha,\beta)$-tour of $\Bc_{{\underline n}}$ exists then $\alpha$
 and $\beta$ must be coprime. Indeed, $(1,1,\ldots,1)$ is connected to
 $(2,1,1\ldots,1)$. Therefore, there are $k,l\in\Z$ such that
 $k\alpha+l\beta=1$.
\end{remark}

We extend the concept of sites to this setting and
  recall that $(e_j)_{j\in \lint 1,   n\rint}$ denotes the canonical
  basis of $\R^n$. 
\begin{definition} 
 Given an $(\alpha,\beta)$-tour $(a^i)_{i\in I}$, we say that it
 contains a well-oriented $\alpha$-site if there are $n,m\in I$ and
$j\in \lint 1,n\rint$, such that
\[a^{n+1}-a^{m}=\pm (a^{m+1}-a^{n}) \in \{-\alpha e_j, \alpha
e_j\}.\] 
We denote it by $(a^n, a^{n+1}, a^{m}, a^{m+1})$. 

Moreoever, we say that it contains a non-well-oriented
  $\alpha$-site if there are $n,m\in I$ and 
$j\in \lint 1,n\rint$, such that
\[a^{n}-a^{m}=\pm (a^{m+1}-a^{n+1}) \in \{-\alpha e_j,  \alpha e_j\}.\]
We denote it by $(a^n, a^{n+1}, a^{m+1}, a^{m})$. We define similarly $\beta$-sites.
\end{definition}

As in Section \ref{s:patterns} we
can use $\alpha$-sites and $\beta$-sites to join together two
cycles. The difference is that it connects two copies which are
  $\beta$ ($\alpha$ resp.) far from one another.

\begin{example}\label{e:alpha}
 More concretely, given an Hamiltonian $(a^i)_{i\in I}$ on
$\Gc^{\alpha,\beta}_{\underline n}$ that contains a well-oriented
$\alpha$-site $(a^n,a^{n+1}, a^m, a^{m+1})$,  let ${\underline
  m}:=({\underline n},\beta + 1)$ so that $\Bc_{\underline m} =
(\Bc_{\underline n},1) \times (\Bc_{\underline n},2) \times \ldots
\times (\Bc_{\underline n},\beta + 1)$. Noting that
$\Gc_{{\underline m }}^{\alpha,\beta}$ contains the two cycles,
$(a^i,1)_{i\in I}$ and  $(a^i,\beta + 1)_{i\in I}$, we see that  
\begin{align*}
(a^{n}, 1)&\to (a^{m+1}, \beta +1 ) \to (a^{m+2}, \beta + 1) \to \ldots \to
(a^{m}, \beta + 1) \to
\\
& \to (a^{n+1}, 1) \to (a^{n+2}, 1)\to \ldots \to (a^{n}, 1)
\end{align*}
is a cycle on $(\Bc_{\underline n},1) \cup (\Bc_{\underline n},\beta +
1)$. The case for a non-well-oriented  site is similar except that we
reverse the orientation of the $(\beta+1)$-th copy.
\end{example}
 We first prove an analogue of Proposition 
\ref{p:gain} for $(\alpha,\beta)$-sites.

\begin{proposition}\label{p:gainab}
Suppose that $\Gc^{\alpha,\beta}_{\underline
  {n}}$ contains a Hamiltonian cycle which contains 
$2$ $\alpha$-sites and $2$ $\beta$-sites, which are two by two
  disjoint. Take now $\underline{m}:=(\underline{n}, k)$, with $k\geq
\alpha+\beta-1$. Then  
$\Gc^{\alpha,\beta}_{\underline{m}}$ contains a Hamiltonian cycle which contains 
$2$ $\alpha$-sites and $2$ $\beta$-sites, which are two by two
  disjoint. 
\end{proposition}

\proof 
Recall that by Remark \ref{r:coprime}  $\alpha$ and $\beta$ are
  coprime. Set $\alpha > \beta$. As the demonstration is
similar, we shall present only the case 
where all $4$ sites are well-oriented. Let $(a^i)_{i\in I}$ be the Hamiltonian
cycle. Note that $(a^i, j)_{i\in I}$, for $j\in \lint 1, k \rint$, are
$k$ disjoint cycles in $\Gc^{\alpha,\beta}_{\underline{m}}$. 

We   start with the two $\alpha$-sites $(a^{n_1}, a^{n_1+1}, a^{m_1}, a^{m_1+1})$ and
$(a^{n_2}, a^{n_2+1},  a^{m_2}, a^{m_2+1})$ in $(a^i)_{i\in I}$. 
As in Proposition
\ref{p:gain} and in the spirit of Example \ref{e:alpha}, for all
$i\in \lint 1, \beta \rint$, we connect the
$i$-th layer to  the $i + \beta$-th layer with the help of the first
$\alpha$-site, then the $i+2\beta$-th layer with the help of the
second $\alpha$-site $\ldots$ until we have a cycle on 
\[\Tc_i:=\bigcup_{c\in \biglint 0, \lfloor \frac{k-i}{\beta} \rfloor
  \beta\bigrint} (\Bc_{{\underline
    n}},\beta c+i),\]  
for all $i\in \lint 1, \beta \rint$. Note that $\cup_{i\in \lint
  1, \beta\rint} \Tc_i$ is a partition of
$\Bc_{\underline{m}}$. More concretely, for the first 
layer and when   $\lfloor \frac{k-1}{\beta} \rfloor$ is even,
this gives:   
\begin{align}
\nonumber
(a^{n_1}, 1)&\to (a^{m_1+1}, \beta + 1) \to (a^{m_1+2}, \beta + 1) \to \ldots \to
(a^{n_2}, \beta + 1) \to
\\
\nonumber
& \to (a^{m_2+1}, 2\beta + 1 ) \to (a^{m_2+2}, 2\beta + 1)\to \ldots \to (a^{n_1}, 2\beta + 1)\to 
\\
\nonumber
& \to (a^{m_1+1}, 3\beta + 1) \to (a^{m_1+2}, 3\beta + 1) \to \ldots \to (a^{n_2}, 3\beta + 1)
\to \ldots
\\
\nonumber
& \to (a^{m_2+1}, \textstyle{(\lfloor \frac{k-1}{\beta} \rfloor -1)} \beta + 1) \to (a^{m_2+2},(\lfloor \frac{k-1}{\beta} \rfloor -1 ) \beta + 1) \to \ldots \to (a^{n_1}, (\lfloor \frac{k-1}{\beta} \rfloor -1 ) \beta + 1) \to
\\
\nonumber
& \to (a^{m_1+1}, \textstyle{\lfloor \frac{k-1}{\beta} \rfloor \beta + 1)} \to (a^{m_1+2}, \lfloor \frac{k-1}{\beta} \rfloor \beta + 1) \to \ldots \to (a^{m_1}, \lfloor \frac{k-1}{\beta} \beta \rfloor + 1)
\to
\\
\nonumber
& \to (a^{n_1+1},\textstyle{(\lfloor \frac{k-1}{\beta} \rfloor -1 )} \beta + 1) \to (a^{n_1+2}, (\lfloor \frac{k-1}{\beta} \rfloor -1 ) \beta + 1) \to \ldots \to (a^{m_2},
(\lfloor \frac{k-1}{\beta} \rfloor -1 ) \beta + 1) \to \ldots
\\
\nonumber
& \to (a^{n_2+1}, \beta + 1) \to (a^{n_2+2}, \beta + 1) \to \ldots \to (a^{m_1}, \beta + 1)
\\
\nonumber
&\to  (a^{n_1+1}, 1)\to (a^{n_1+2}, 1)\to \ldots \to (a^{n_1}, 1)
\end{align}
Note that $\big((a^{n_2},1),
(a^{n_2+1},1), (a^{m_2},1), (a^{m_2+1},1)\big)$ and
$\big((a^{n_2},\lfloor \frac{k-1}{\beta} \rfloor \beta + 1),
(a^{n_2+1},\lfloor \frac{k-1}{\beta} \rfloor \beta + 1),
\linebreak
(a^{m_2},\lfloor \frac{k-1}{\beta} \rfloor \beta + 1),
(a^{m_2+1},\lfloor \frac{k-1}{\beta} \rfloor \beta + 1)\big)$ are two ``free''
 $\alpha$-sites. When $\lfloor
 \frac{k-1}{\beta} \rfloor$ is odd, we replace the latter by
$\big((a^{n_1},\lfloor
\frac{k-1}{\beta} \rfloor \beta + 1), 
(a^{n_1+1},\lfloor \frac{k-1}{\beta} \rfloor \beta + 1),
(a^{m_1},\lfloor \frac{k-1}{\beta} \rfloor \beta + 1),
(a^{m_1+1},\lfloor \frac{k-1}{\beta} \rfloor \beta + 1)\big)$. 

We now use the two $\beta$-sites $(a^{p_1}, a^{p_1+1}, a^{q_1},
a^{q_1+1})$ and $(a^{p_2}, a^{p_2+1}, a^{q_2}, a^{q_2+1})$ in
$(a^i)_{i\in I}$ to join $\Tc_1,\ldots, \Tc_{\beta}$ into
a Hamlitonian cycle. To lighten notation we denote by $[p_1, i]$ (resp.\
$[p_2, i]$) for $i\in \lint 1, k\rint$, the $i$-th copy of the first
(resp.\ the second) $\beta$-site. Let $d_1 := 1$ and, for $i\in \lint
2, \beta\rint$ we set $d_i\in \lint 1, \beta\rint$ such that $d_{i} =
d_{i-1} + \alpha \mod{\beta}$. Since $\alpha$ and $\beta$ are coprime,
we stress that $d_i$  is well-defined and that the map $i\mapsto d_i$
is a bijection onto $\lint 1, \beta\rint$. 

Using $[p_1, d_1]=[p_1, 1]$ with $[p_1, d_1+\alpha]$, we connect
$(\Bc_{\underline{n}},1)$ to $(\Bc_{\underline{n}},1+\alpha)$.   We
have then constructed a Hamiltonian cycle for 
$\Tc_{d_1}\cup\Tc_{d_2}$. The
$\beta$-sites $[p_2,d_1]$ and $[p_2, d_1+\alpha]$ will not be used
anymore. Inductively, we connect $\Tc_{i-1}$ to
$\Tc_i$  for all $i\in \lint 2, \beta\rint$ using $[p_1, d_{i-1}]$
with $[p_1, d_{i-1}+\alpha]$.  Note that it is possible since $k\geq
\alpha+\beta-1$ and since $d_i\neq d_{i-1}+\alpha$, recall that
$\alpha>\beta$. We conclude by recalling that $\cup_{i\in   \lint 1,
  \beta\rint} \Tc_i = \cup_{i\in   \lint 1, \beta\rint} \Tc_{d_i}$. \qed

We turn to the main result of this section. As in the
$(1,2)$-tour case, it is sufficient to construct tours in a low
dimension to show the existence of tours on all, sufficiently large,
chessboards in higher dimensions.  
\begin{theorem}\label{t:ab}
Given $\alpha > \beta$, $n_1,n_2 \geq 2\alpha +1$ and $n_3,\ldots, n_k
\geq \alpha + \beta -1$, if an $(\alpha,\beta)$-tour on an $n_1 \times
n_2$ chessboard exists then an $(\alpha,\beta)$-tour on an $n_1 \times
n_2\times \ldots \times n_k$ chessboard exists. 
\end{theorem}

\proof We proceed by induction on: 
\\
\\
$\Pc_k:=$``Given a multi-index $\underline{n}=(n_1, n_2, \ldots,
n_{k})$, such that  
\begin{enumerate}
\item $\Gc^{\alpha,\beta}_{(n_1,n_2)}$ contains a Hamiltonian cycle,
\item $n_1,n_2 \geq 2 \alpha +1$,
\item $n_i \geq \alpha+\beta-1$ for all $i \in \lint 3, k \rint$,
\end{enumerate}
 then, $\Gc^{\alpha,\beta}_{\underline{n}}$ contains a Hamiltonian
 cycle which has at least $2$ $\alpha$-sites and $2$ $\beta$-sites." 
\\
\\
for all $k\geq 2$.
\\

\noindent{\bf Basis:} It is a straightforward generalization of
Proposition \ref{p:2dbi} that for $n_1,n_2 \geq 2 \alpha +1$ every
$(\alpha,\beta)$-tour on an $n_1 \times n_2$ chessboard will contain
$2$ $\alpha$-sites and $2$ $\beta$-sites, hence $\Pc_2$ holds.   
\\

\noindent{\bf Inductive step:} Suppose that $\Pc_k$ holds true for some $k\geq
2$. Take $\underline n=(n_1, n_2, \ldots, n_{\bf k+1})$ such that that a), b) an c) hold true. We set 
\begin{align*}
\underline{n'}= (n_1,n_2 \ldots , n_k).  
\end{align*}
We see that ${\underline n'}$ satisfies the hypotheses of $\Pc_k$ and
therefore, by $\Pc_k$, we see that
$\Gc^{\alpha,\beta}_{\underline{n'}}$ contains a Hamiltonian cycle
which has at least $2$ $\alpha$-sites and $2$ $\beta$-sites. . We
apply Proposition \ref{p:gainab} and get that $\Pc_{k+1}$ is true.  

Therefore, we have proved by induction that $\Pc_k$ is true for all
$k\geq 2$. In particular, this proves the theorem. \qed

It is not known in general for which $\alpha,\beta$
$(\alpha,\beta)$-tours exist on 
sufficiently large chessboards. In light of the conditional nature of
Theorem \ref{t:ab} it seems natural to conjecture

\begin{conjecture}
Take $\alpha,\beta$ such that
$\gcd(\alpha+\beta,\alpha-\beta)=1$. Then,   
there exists $M$ such that
an $(\alpha,\beta)$-tour exists on all $n_1 \times n_2$ chessboards,
where $n_1$ is even and $n_1,n_2 \geq M$.
\end{conjecture}

\appendix 
\renewcommand{\theequation}{\thesection .\arabic{equation}} 

\section{An alternative construction of Tours in dimension $3$}
\subsection{Around the construction of Schwenk}
For the sake of completeness, we discuss the Theorem of Schwenk.

\begin{theorem}\label{t:S2}
 Let $1\leq m\leq n$. The $m\times n$
 bi-sited chessboard has a closed knight tour if and only if
\begin{enumerate}
\item $m$ or $n$ is even,
\item $m\notin\{1,2,4\}$,
\item $(m,n)\neq (3,4), (3,6)$ or $(3,8)$. 
\end{enumerate} 
\end{theorem} 
At the end of the section, we give a proof for the existence of
closed tours in the setting of Theorem \ref{t:S}. We refer to
\cite{Sch, Wat} for the proof of the non-existence of tours.
Note that by Proposition \ref{p:2dbi}, their structure would be used to
reduce the number of cases to study in dimension $3$. We start with
two definitions.  


\begin{definition} Given a board of size $(m,n)$, a (open or closed)
  tour is called \emph{seeded} if it 
  includes the edges $((1,m-2),(2,m))$ and $((n-2,1),(n,2))$. 
\end{definition}

Some examples are given in Figures \ref{f:3x10} and \ref{f:3x12} below. 

\begin{remark}\label{r:seeded}
A board of size $(m,n)$ has an open (resp.\ closed) seeded tour if and
only if one of size $(n,m)$ has.  
\end{remark} 

\begin{definition}
Given a board of size $(4,m)$, an open tour is called a $4 \times m$ 
\emph{extender} if the tour  starts at $(4,m)$ and ends at
$(4,m-1)$.
\end{definition}

We start by showing the existence. 

\begin{lemma}\label{l:extender}
There exists a seeded $4 \times m$ extender for all  $m \neq 1,2$ or
$4$.
\end{lemma} 

\begin{figure}[ht]
\centering
\includegraphics[scale=0.2]{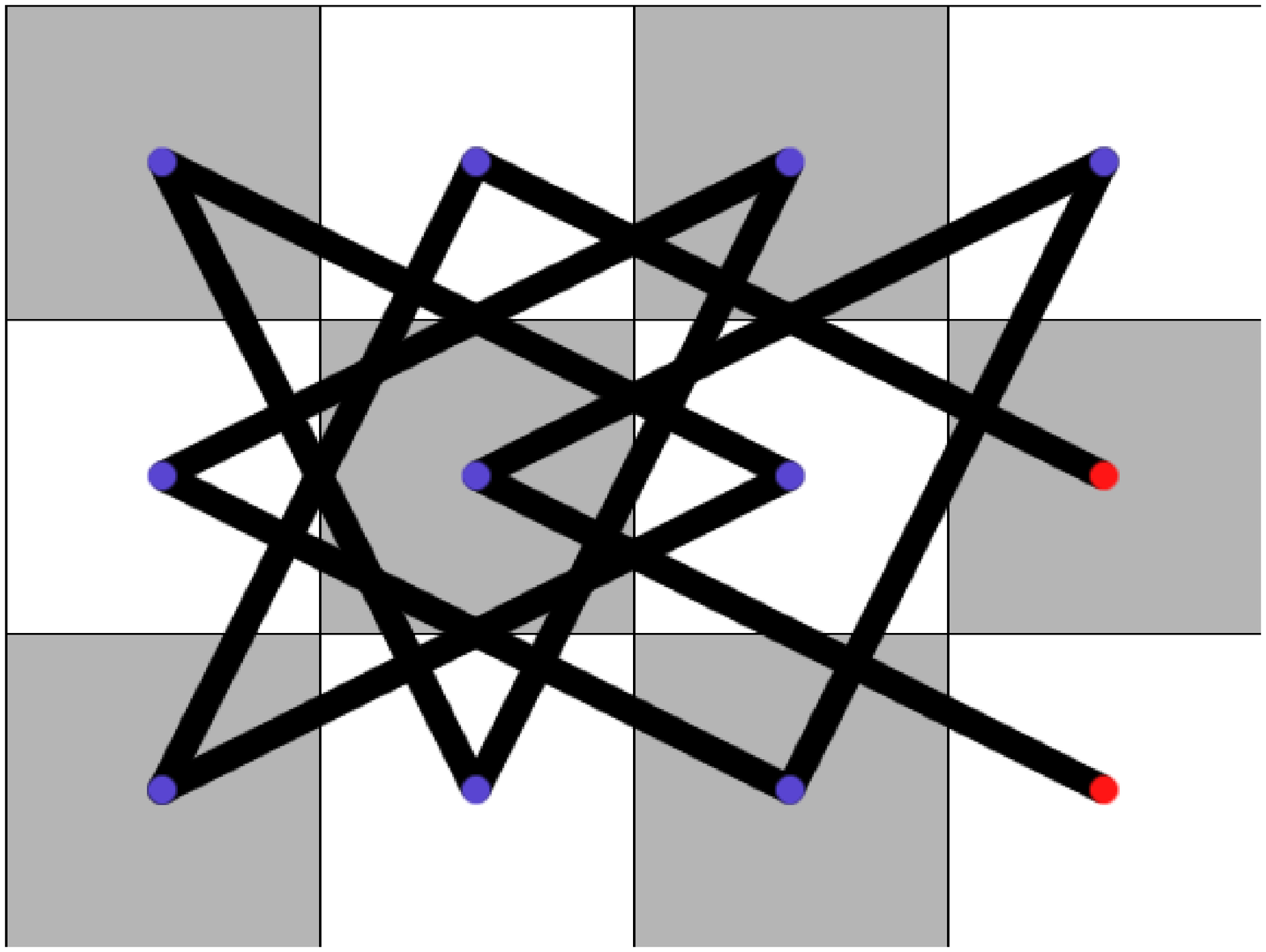}\quad \includegraphics[scale=0.2]{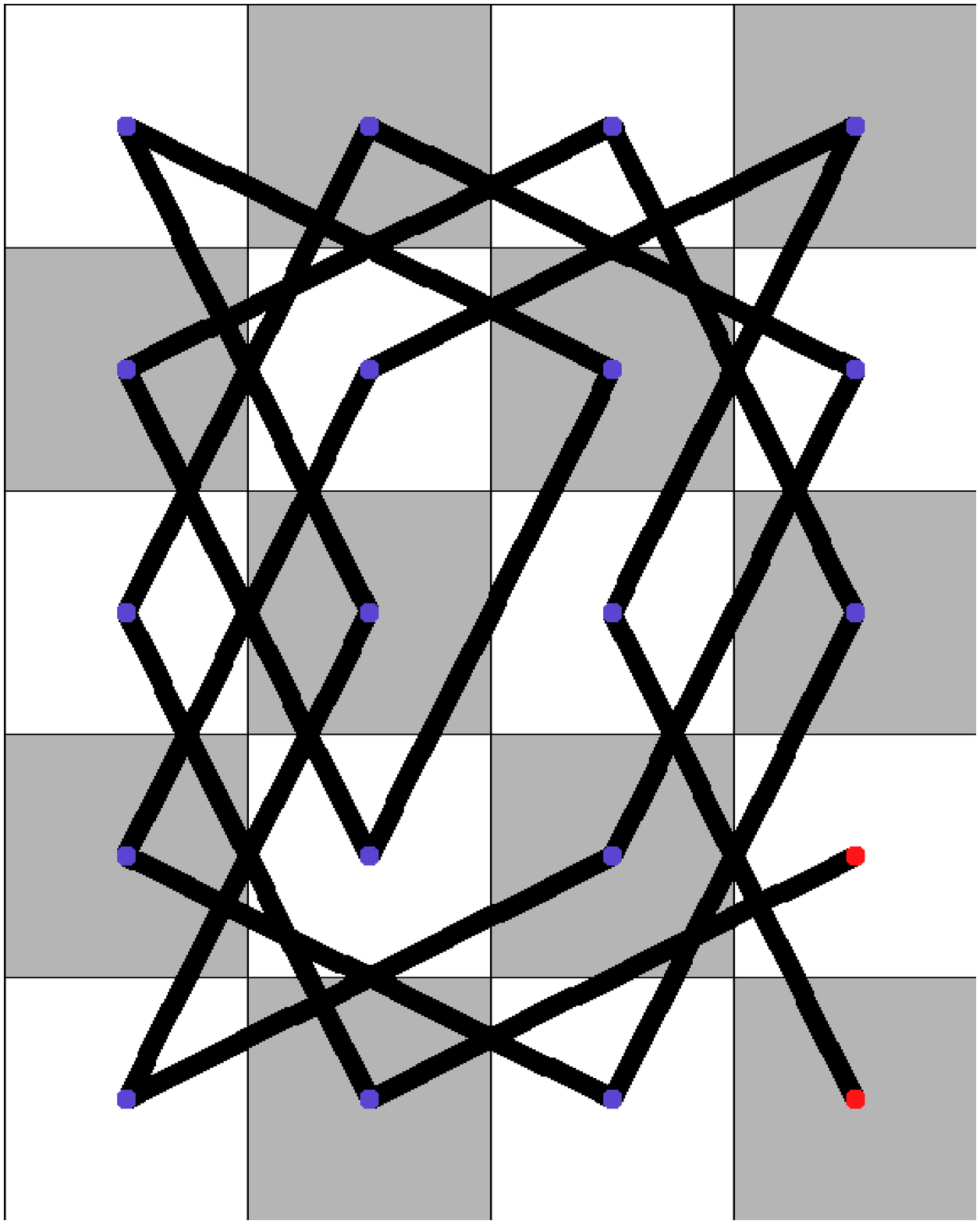}
\quad \includegraphics[scale=0.2]{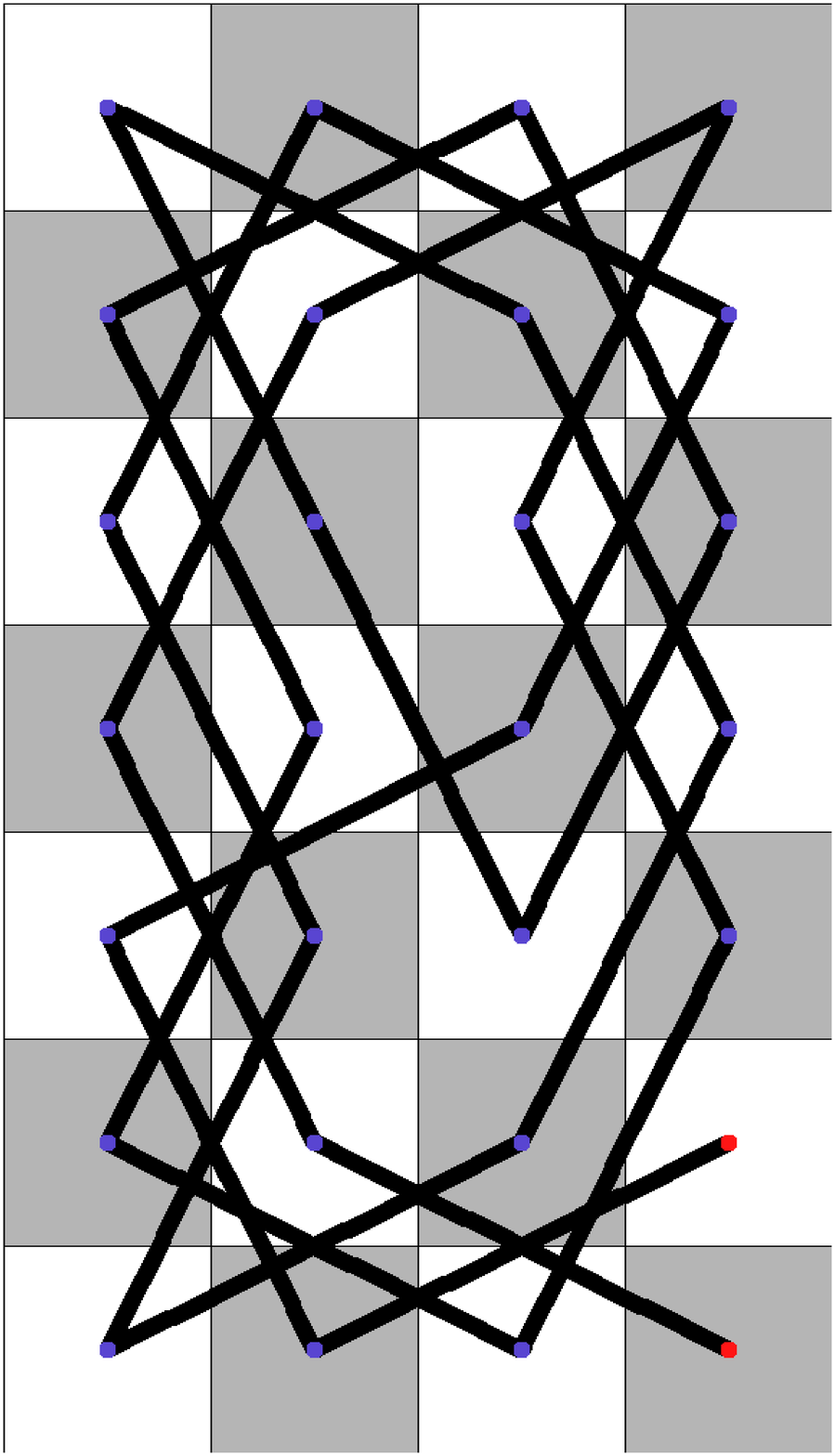} 
\caption{A $4 \times 3$, $4 \times 5$, and a $4 \times 7$ extender.}\label{f:3}
\end{figure}

\begin{figure}[ht]
\centering
\includegraphics[scale=0.2]{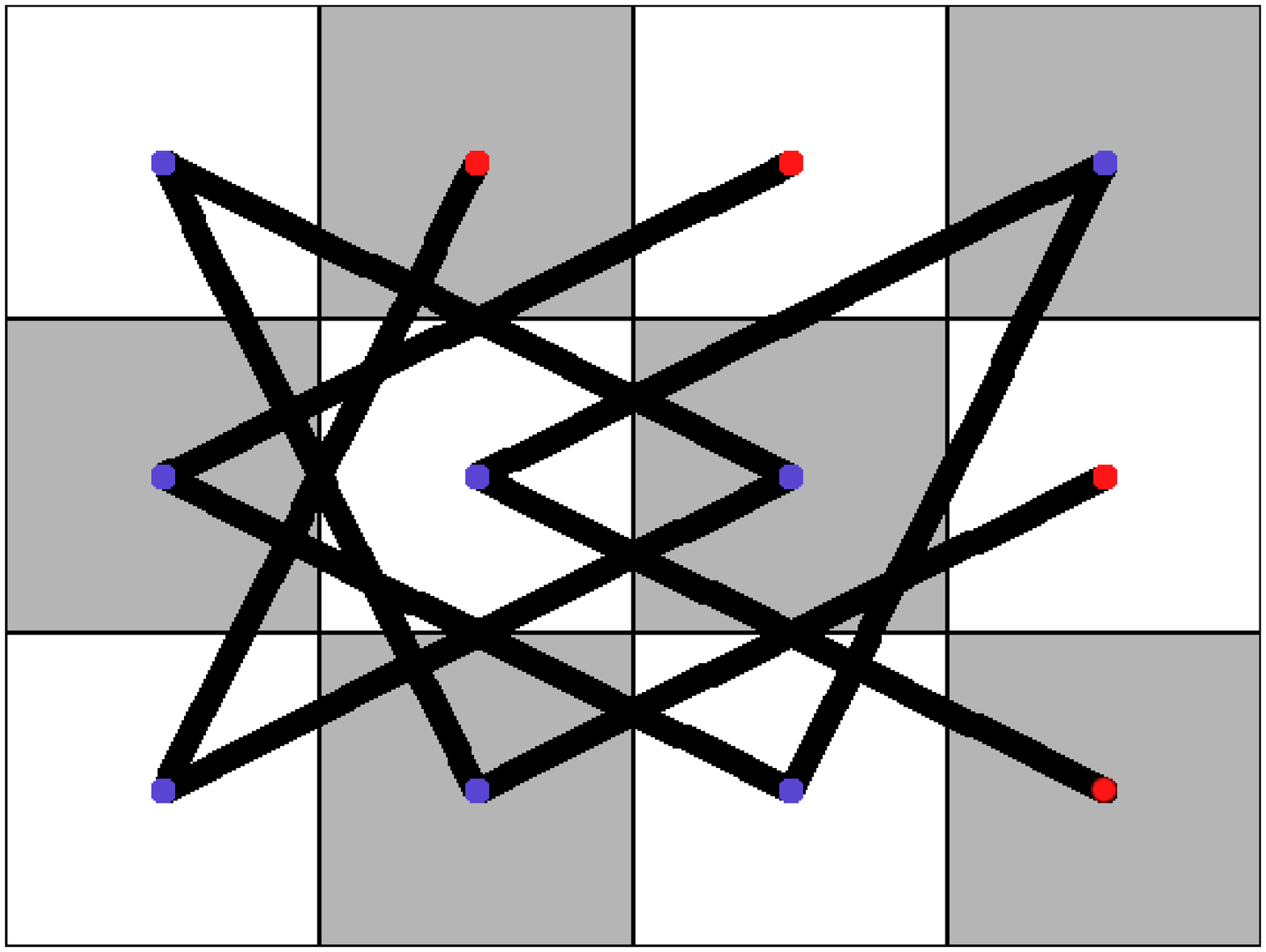}
\caption{An extending pattern.}\label{f:ext}
\end{figure}

\proof The three elementary cases are given in Figure \ref{f:3}.
Then observe that if we place the Figure \ref{f:ext} below a seeded $4
\times m$ extender and add the lines $((4,m-1), (3,m+1))$ and $((4,m),
(2,m+1))$ then it  will form a seeded $4 \times (m+3)$
extender. Conclude by induction. 
\qed

\begin{proposition}\label{p:ext}
If a seeded  $n \times m$ closed tour exists then a seeded $(n + 4k)
\times (m+4l)$ closed tour exists, for all $k,l\in \N$.
\end{proposition}

\begin{figure}[ht]
\centering
\includegraphics[scale=0.2]{3x4.eps} \quad \includegraphics[scale=0.2]{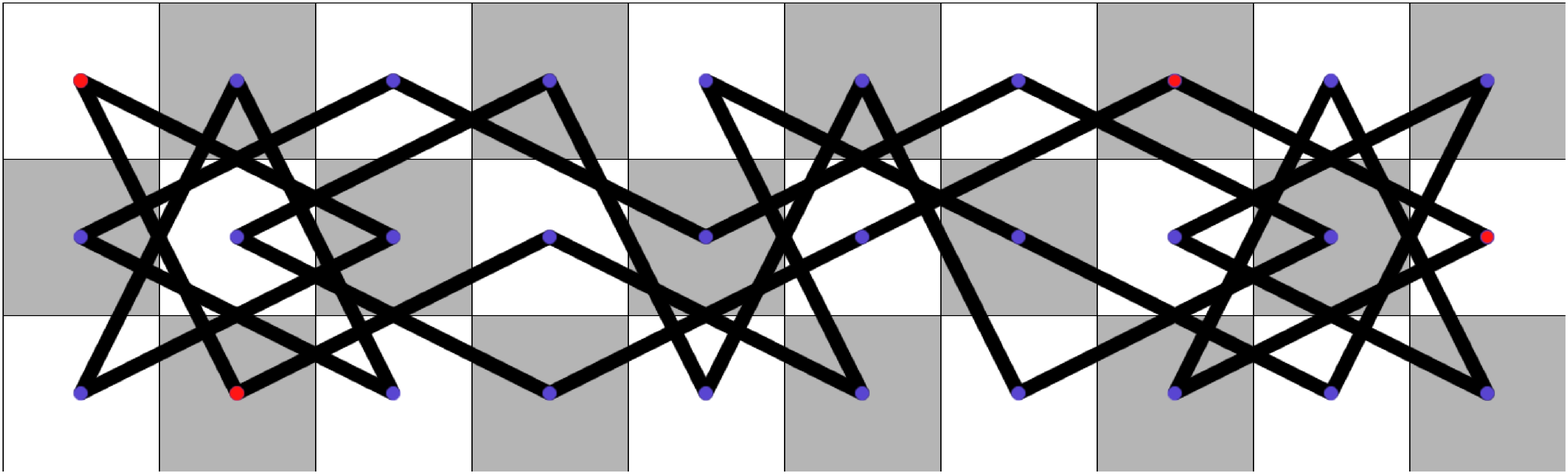}
\caption{How to extend a closed seeded tour}\label{f:extsee}
\end{figure}

\proof
By Theorem \ref{t:S}, we get $m,n \neq 1,2$
or $4$.  Then, there
exists a seeded $4 \times m$ extender. Now if we place a seeded $4
\times m$ extender to the left of a seeded $n \times m$ tour, as in
Figure \ref{f:extsee}.
By removing the line $( (6, m),(5,m-2))$ and adding in the
two lines $((4,m),(5,m-2))$ and $((6,m),(4,m-1))$ we form a $(n + 4)
\times m$ tour. Note also 
that this tour is still seeded. Hence a seeded $(n + 4) \times m$ tour
exists. By induction we get a seeded $(n + 4k) \times m$ tour. Then,
by using Remark \ref{r:seeded} and by repeating the proof, we conclude. \qed 

We are now ready to prove the result.

\begin{figure}[ht]
\centering
\includegraphics[scale=0.2]{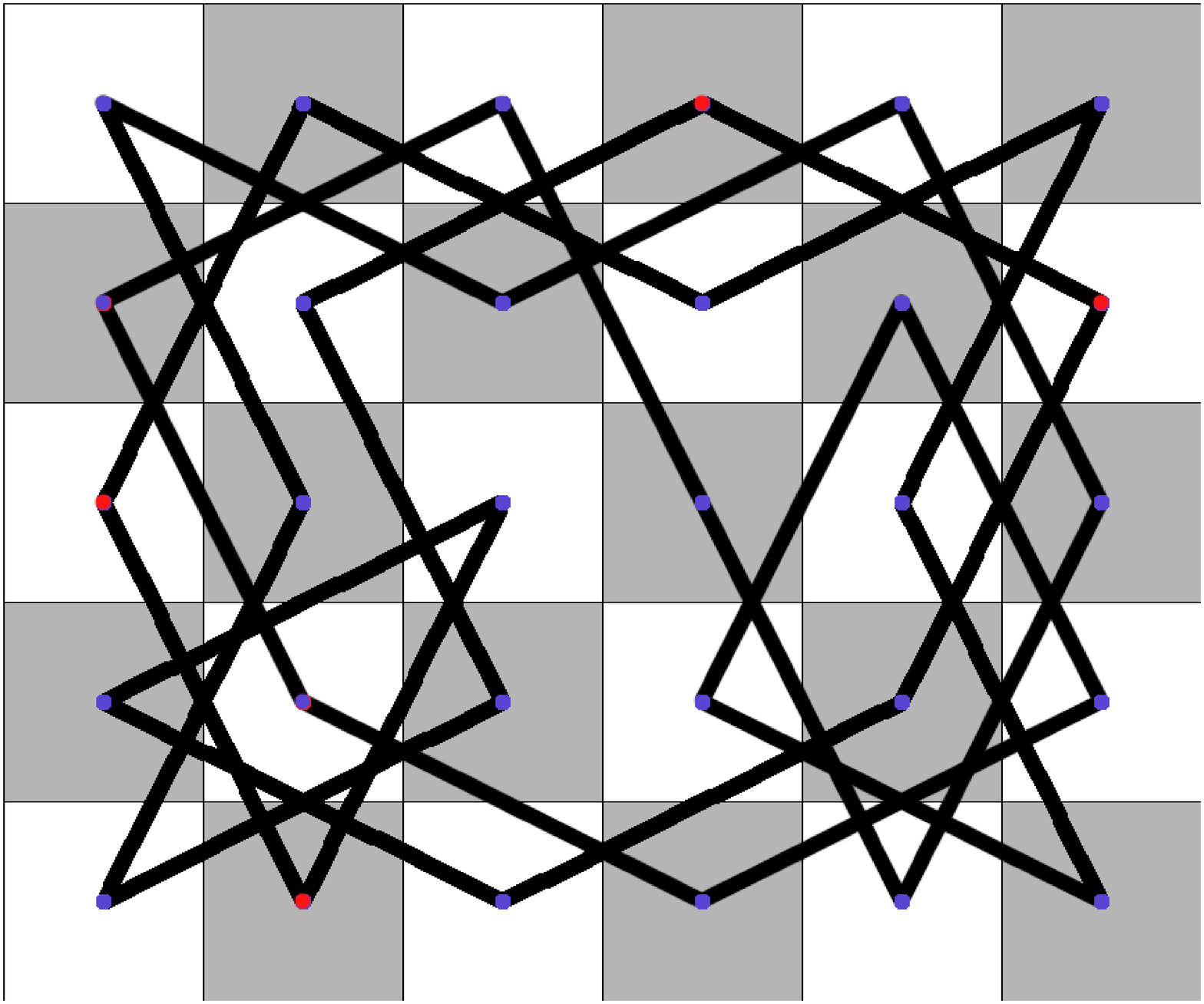}
\quad
\includegraphics[scale=0.2]{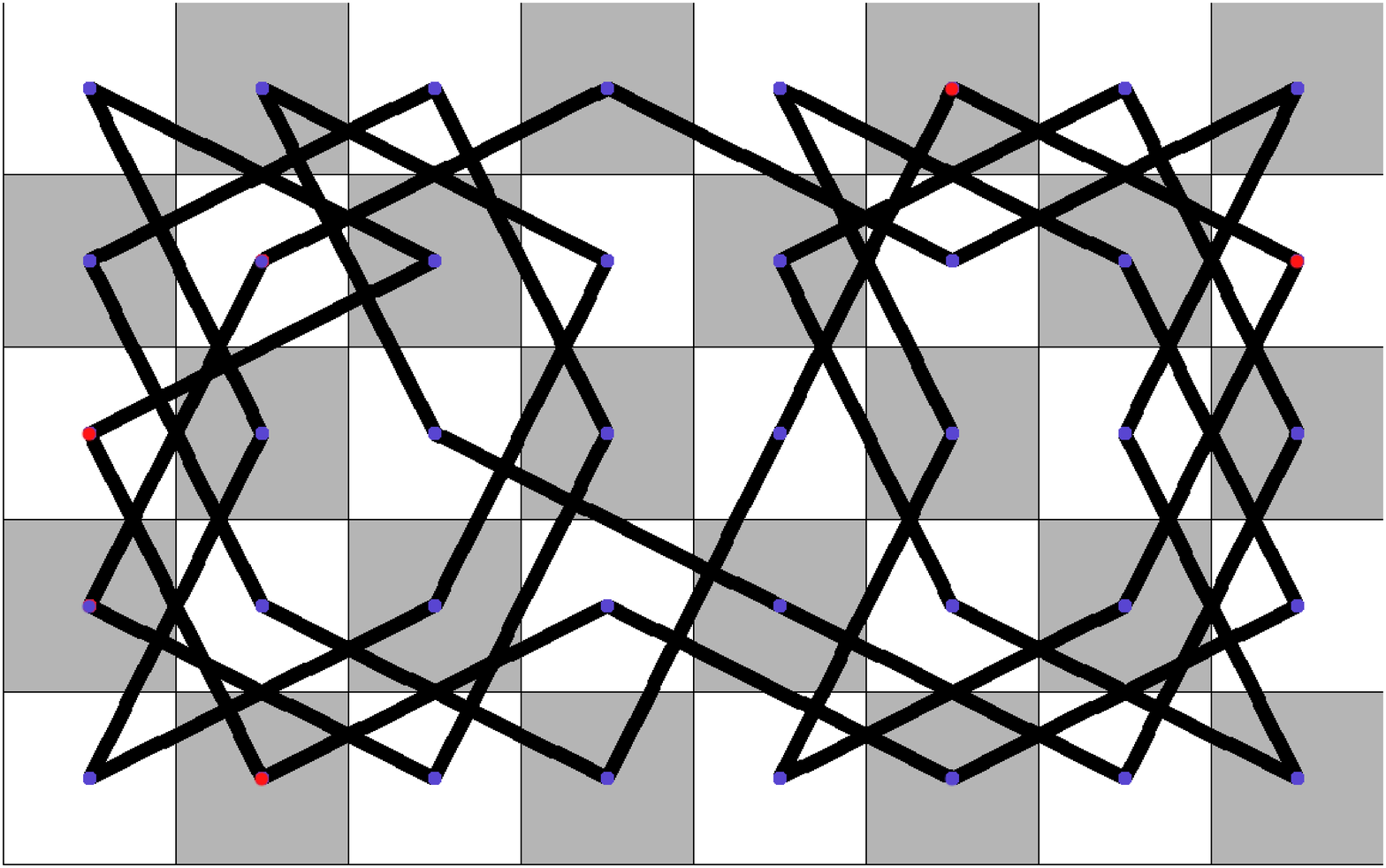}
\caption{A $5\times 6$ and a $5\times 8$ seeded closed tour}\label{f:5x6-5x8}
\end{figure}

\begin{figure}[ht]
\centering
\includegraphics[scale=0.2]{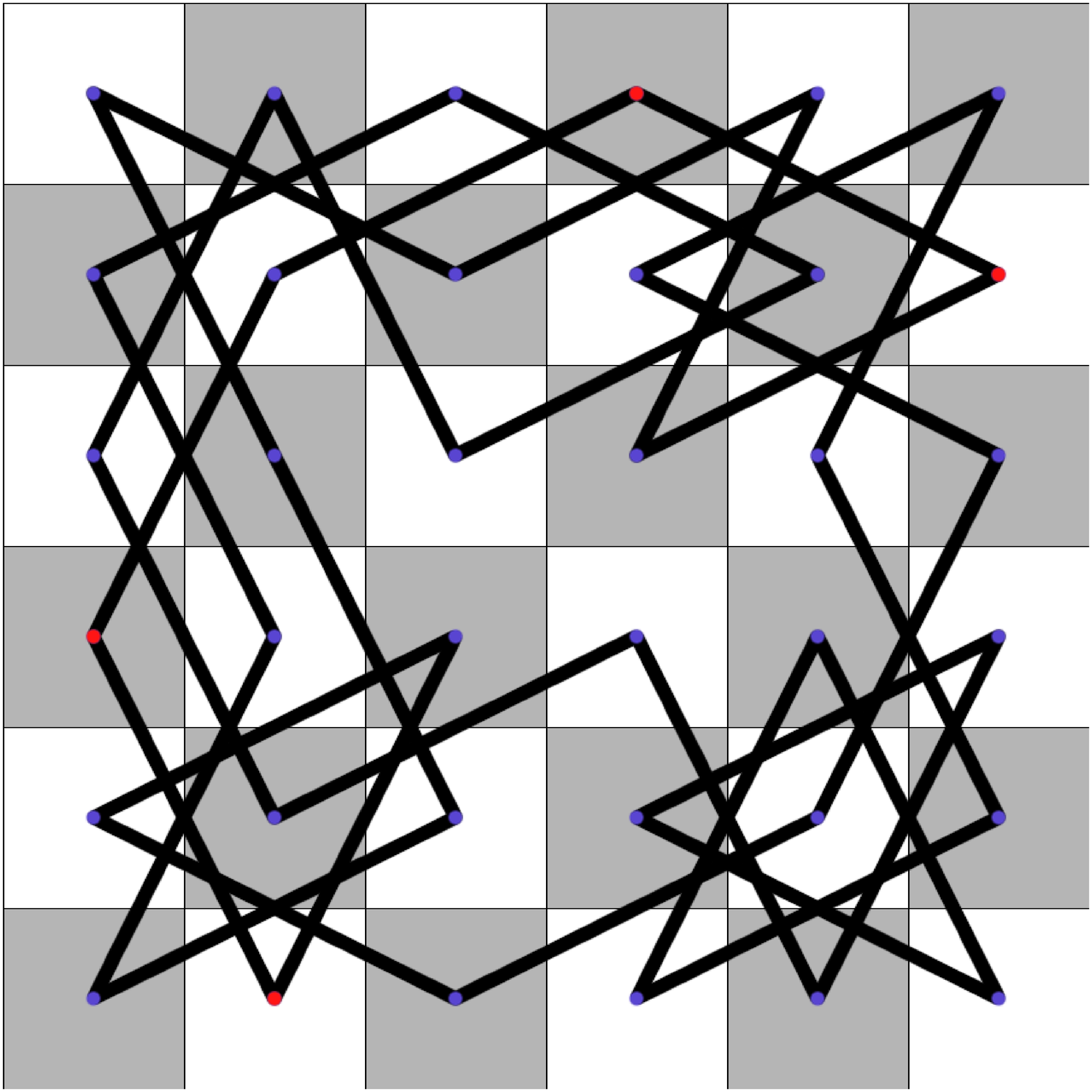}
\quad
\includegraphics[scale=0.2]{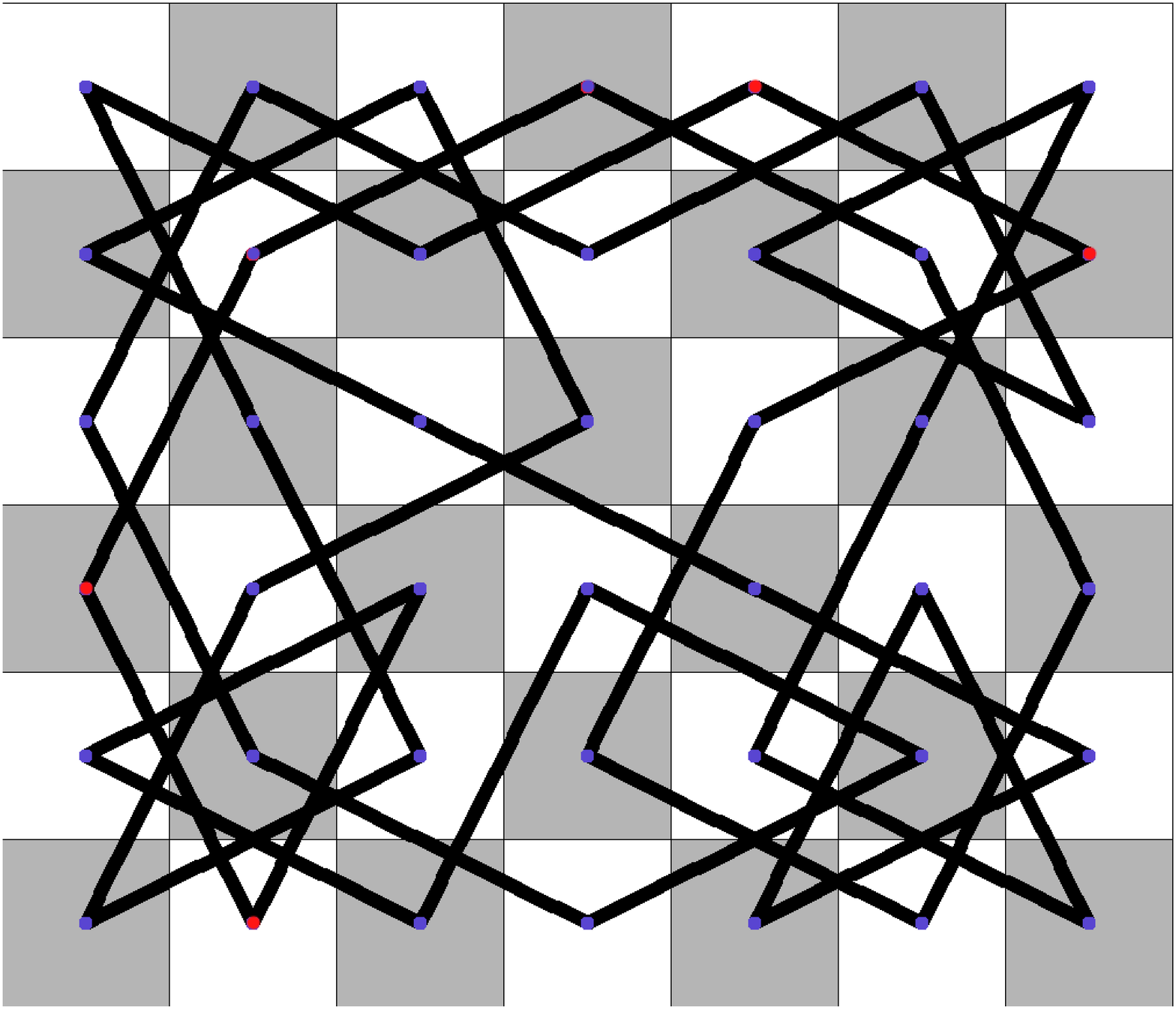}
\caption{A $6\times 6$ and a $6\times 7$ seeded closed tour}\label{f:6x6-6x7}
\end{figure}

\proof[Proof of Theorem \ref{t:S2}]
By Proposition \ref{p:ext}, it is sufficient to exhibit a seeded $n
\times m$ tour for all different pairs of residue modulo $4$
(excepting the cases where both are odd), and possibly some small
cases. A quick check will show it is enough to use as base cases
seeded $3 \times 10$, $3 \times 12$, $5 \times 6$, $5 \times 8$, $6
\times 6$, $6 \times 7$, $6 \times 8$, $7 \times 8$ and $8 \times 8$
tours, which appear in Figures \ref{f:3x10}--\ref{f:6x6-6x7}.
\ \\
\begin{figure}[ht]
\centering
\includegraphics[scale=0.2]{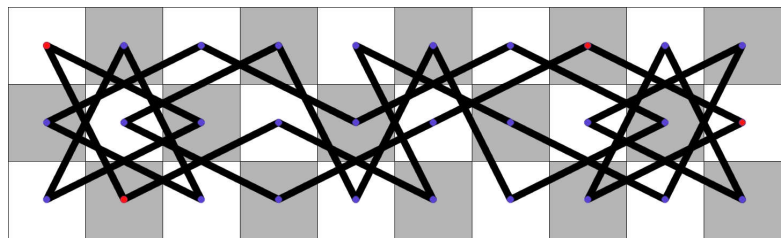}
\caption{A $3 \times 10$ seeded tour.}\label{f:3x10}
\end{figure}

\begin{figure}[ht]
\centering
\includegraphics[scale=0.2]{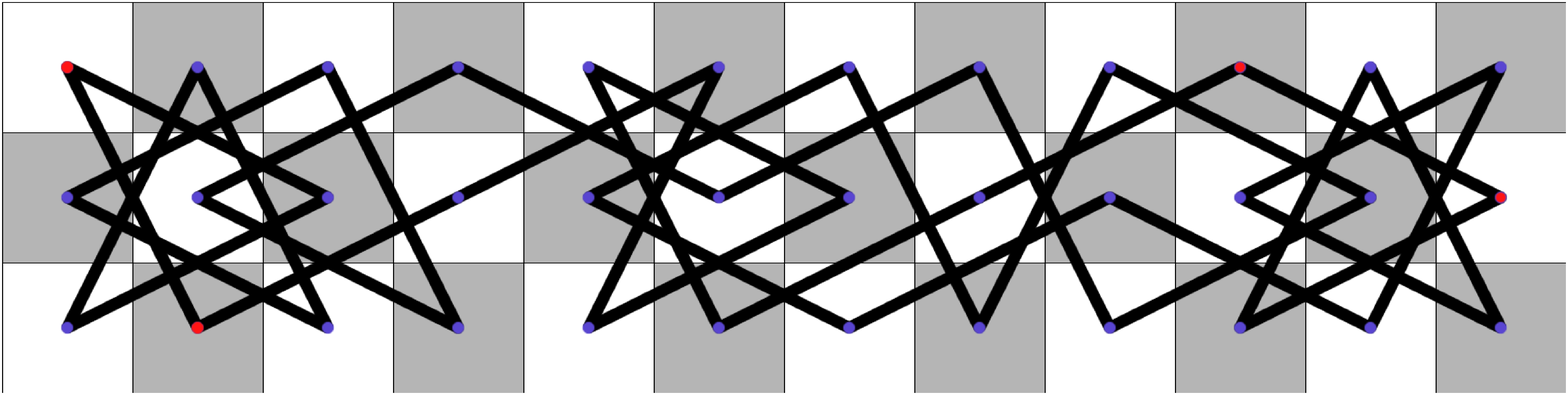}
\caption{A $3\times 12$ seeded closed tour}\label{f:3x12}
\end{figure}

\begin{figure}[ht]
\centering
\includegraphics[scale=0.18]{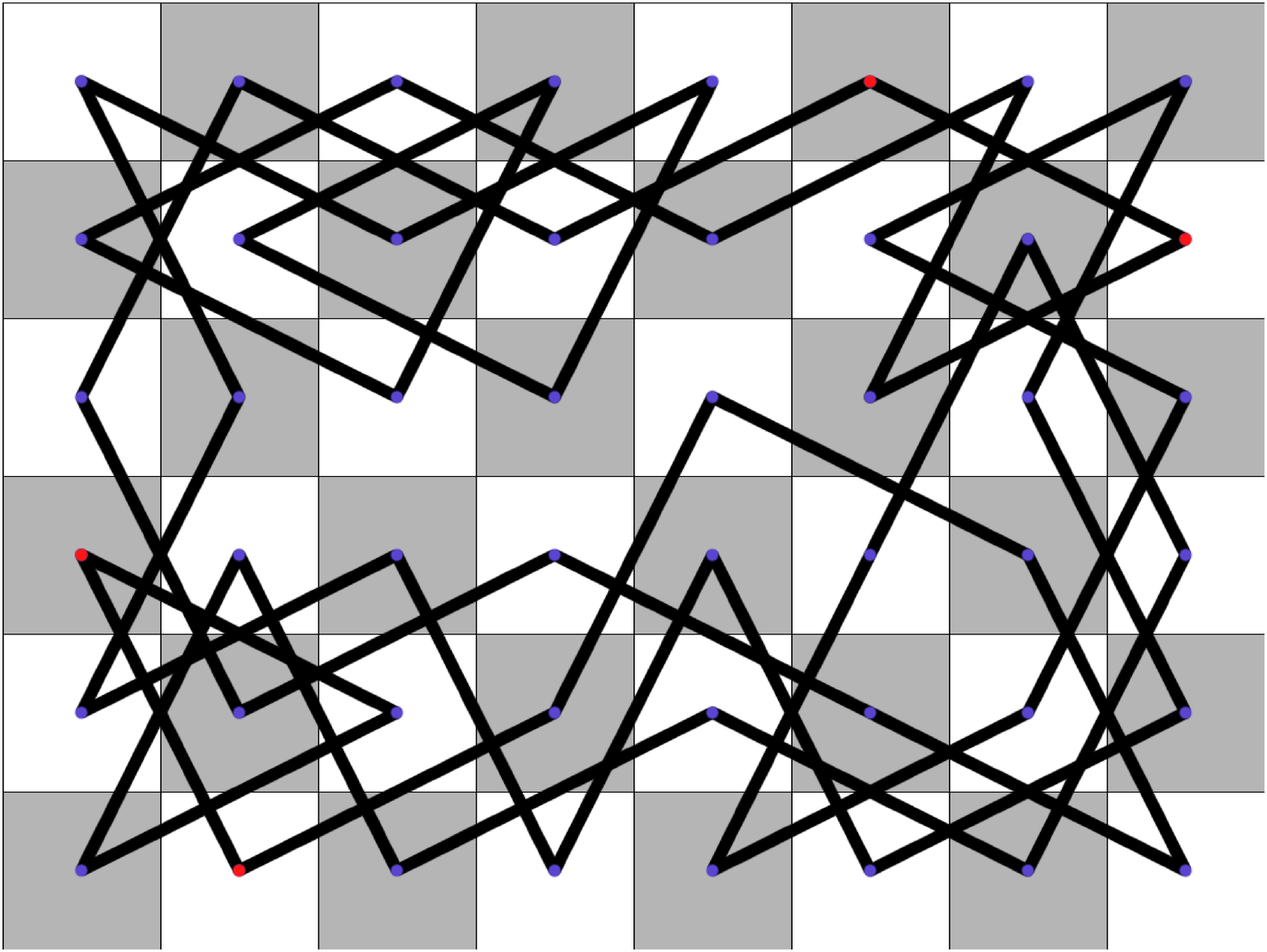}
\quad \includegraphics[scale=0.18]{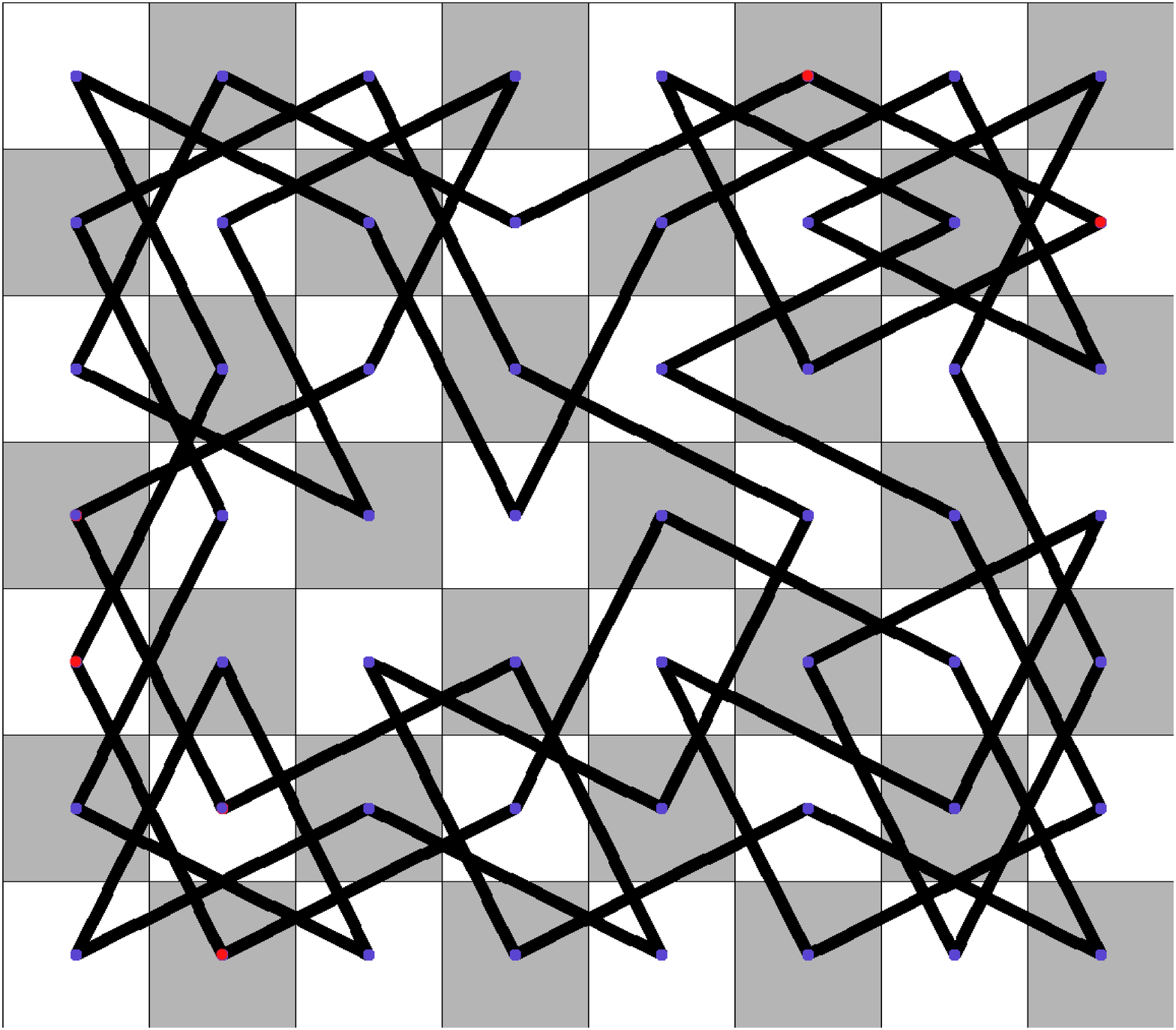} 
\caption{A $6\times 8$ and a $7\times 8$ seeded closed tour}
\end{figure}

\begin{figure}[ht]
\centering
\includegraphics[scale=0.18]{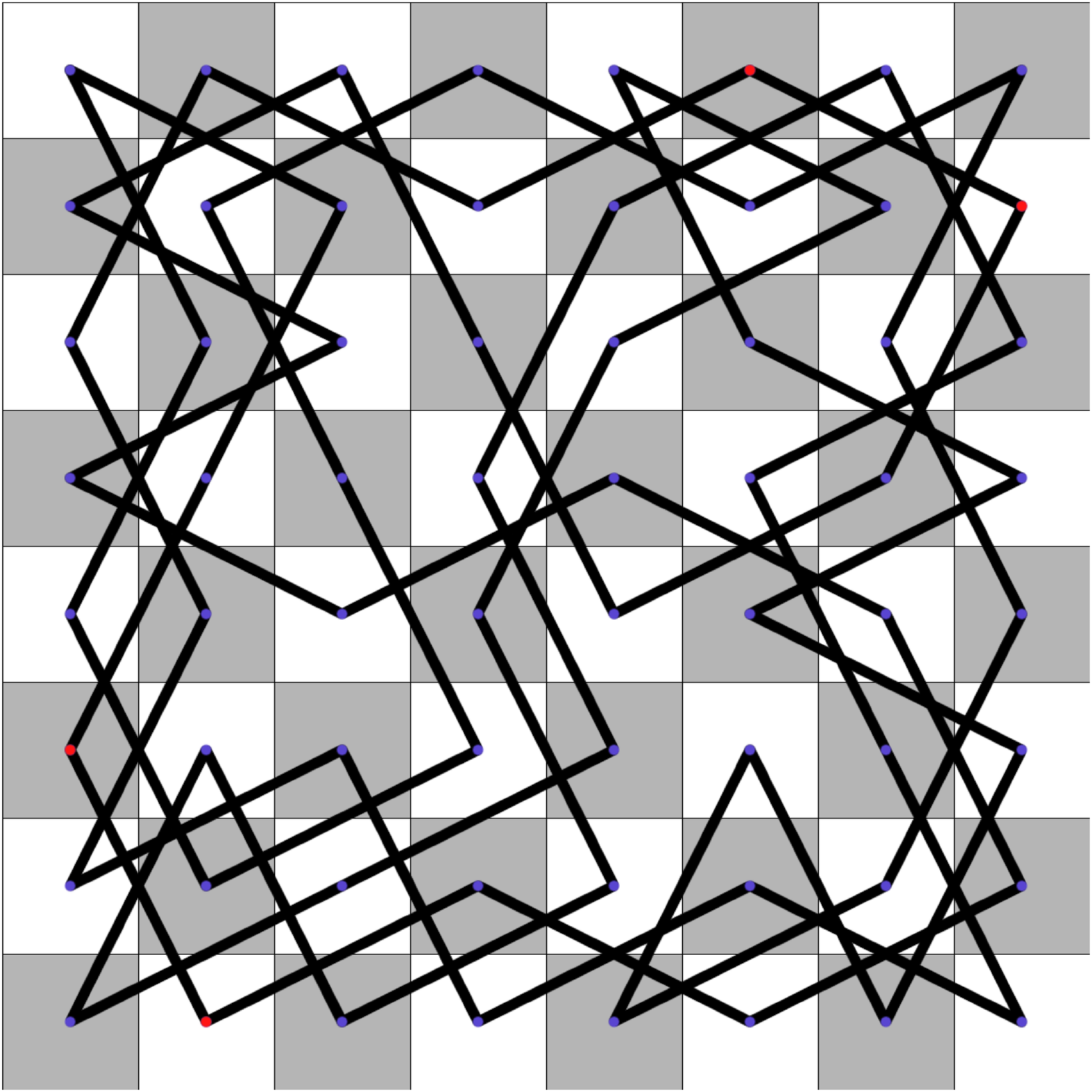}
\caption{A $8\times 8$ seeded closed tour}\label{f:8x8}
\end{figure}

To conclude, recall that Proposition \ref{p:2dbi} ensures that the tours are
all bi-sited.
\qed

\subsection{An altenative proof for $3$-dimensional
  tours}\label{s:DeMalt} 
We are now in position to give an alternative proof of the result of
DeMaio-Mathew and more precisely of Theorem \ref{t:DeM2}. To reduce the
number of cases we have to consider, we combine Theorem \ref{t:S2} and
Proposition \ref{p:gain}.  The non-tourable boards are discussed in Section
\ref{s:forbidden}. It remains to discuss the existence. In the remaining figures 
we will indicate the sites in the tours in red.

\begin{figure}[ht]
\centering
\includegraphics[scale=0.2]{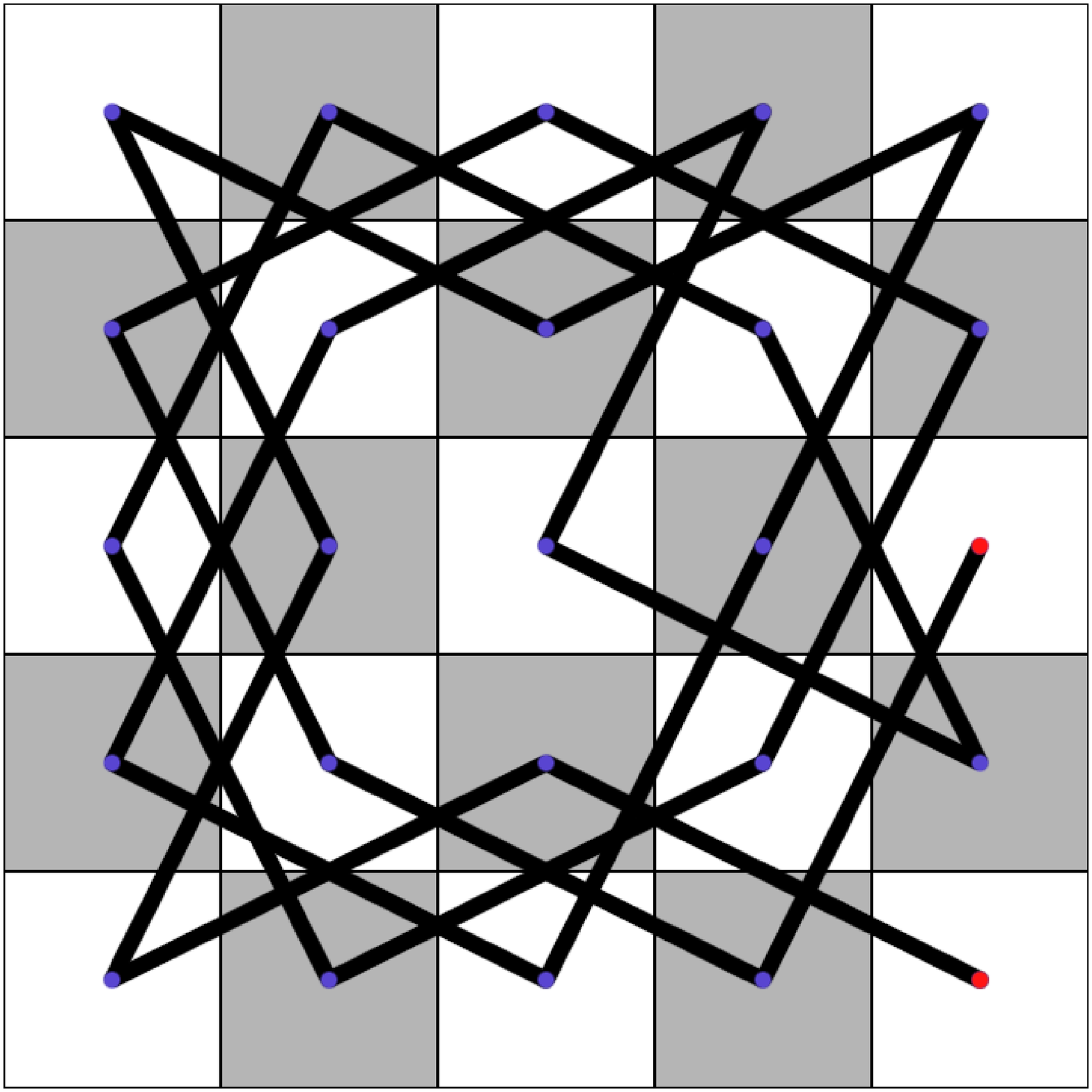}
\quad \includegraphics[scale=0.2]{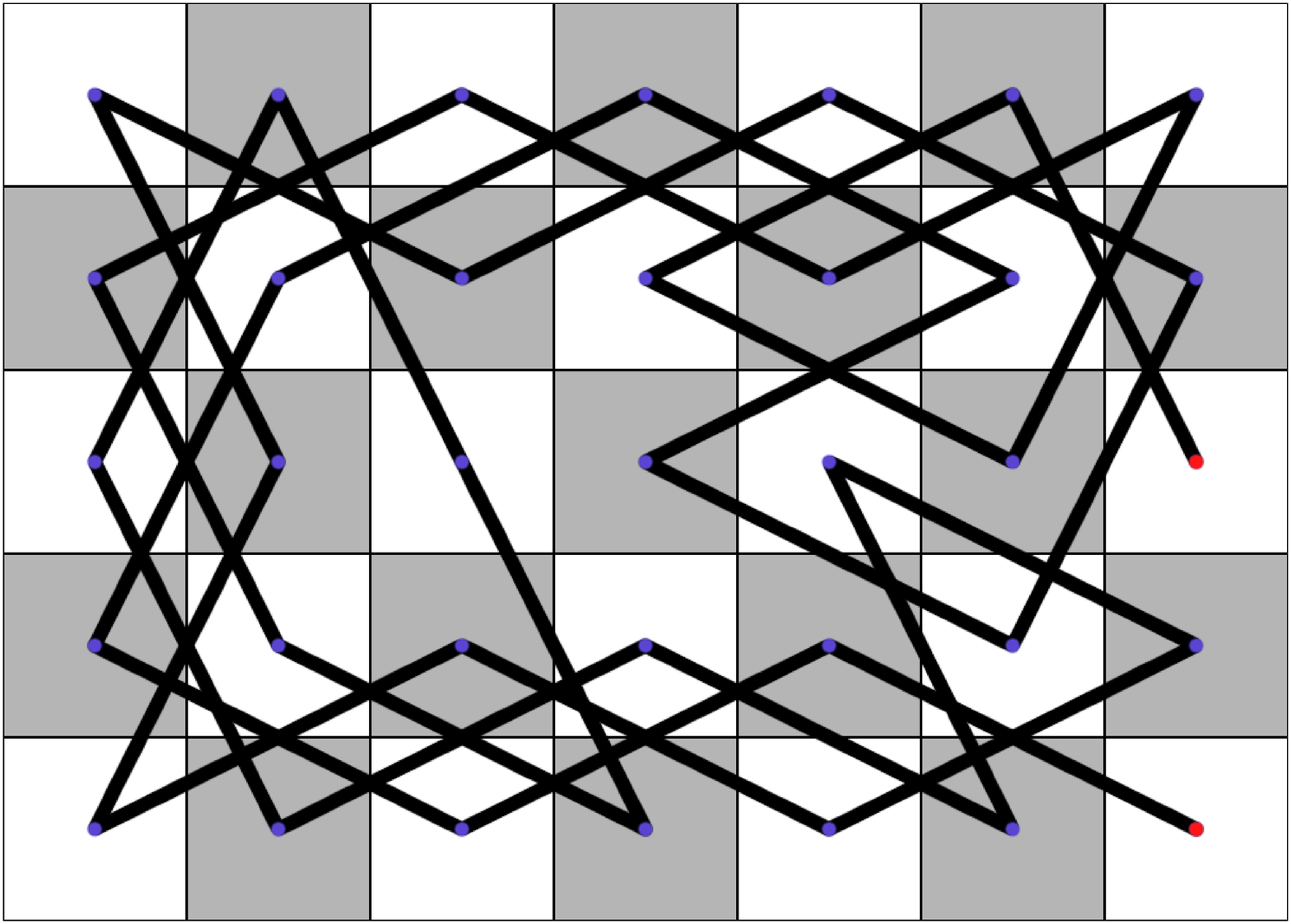} 
\caption{A $5\times 5$ and a $5\times 7$ seeded and bi-sited open
  tour}\label{f:5x5x2and5x7x2} 
\end{figure}

\begin{figure}[ht]
\centering
\includegraphics[scale=0.2]{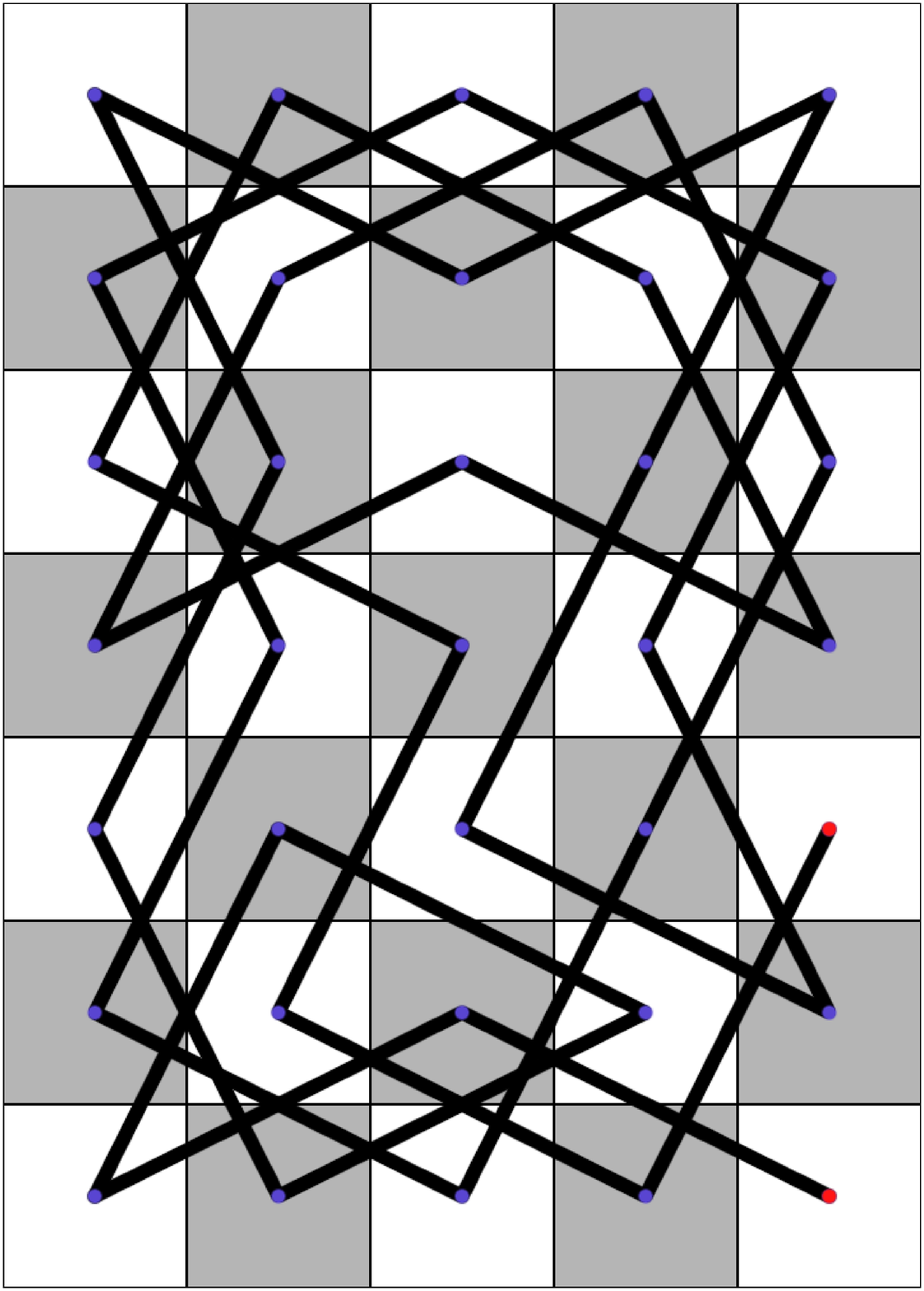}
\quad 
\includegraphics[scale=0.2]{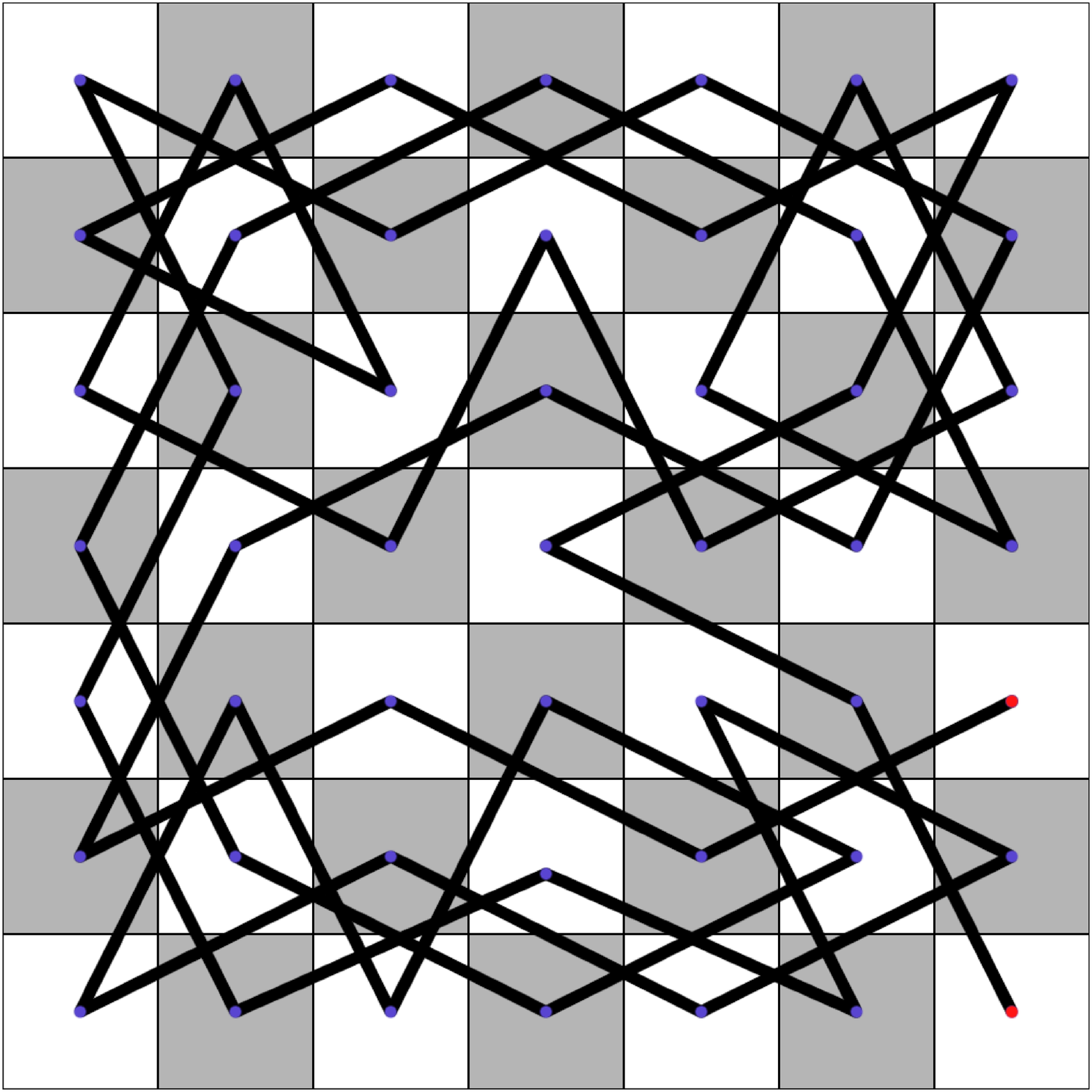}
\caption{A $7\times 5$ and a $7\times 7$ seeded and bi-sited open
  tour}\label{f:7x5x2and7x7x2} 
\end{figure}

\proof
For all $p\in \N\setminus\{0\}$, Theorem \ref{t:S2} and Proposition
\ref{p:gain} ensure the existence of $3$-dimensional bi-sited tours
for all $n \times m \times p$ when an $n \times m$ chessboard admits a
knight's tour. Moreover, using Remark \ref{r:permu}, this also holds
true for any permutation of $n,m$ and $p$. We will split the remaining
tours into cases. 
\\
\\
{\bf i) Case $\bf n \times m \times 2k$, for $\bf n,m \geq 5$ and odd:}
We start by seeded and bi-sited open tours of size $n \times m$, with
$n,m\in \{5,7\}$ that starts at $(n,m)$ and ends two squares above at
$(n,m-2)$, see Figures \ref{f:5x5x2and5x7x2} and \ref{f:7x5x2and7x7x2}.

\begin{figure}[ht]
\centering
\includegraphics[scale=0.2]{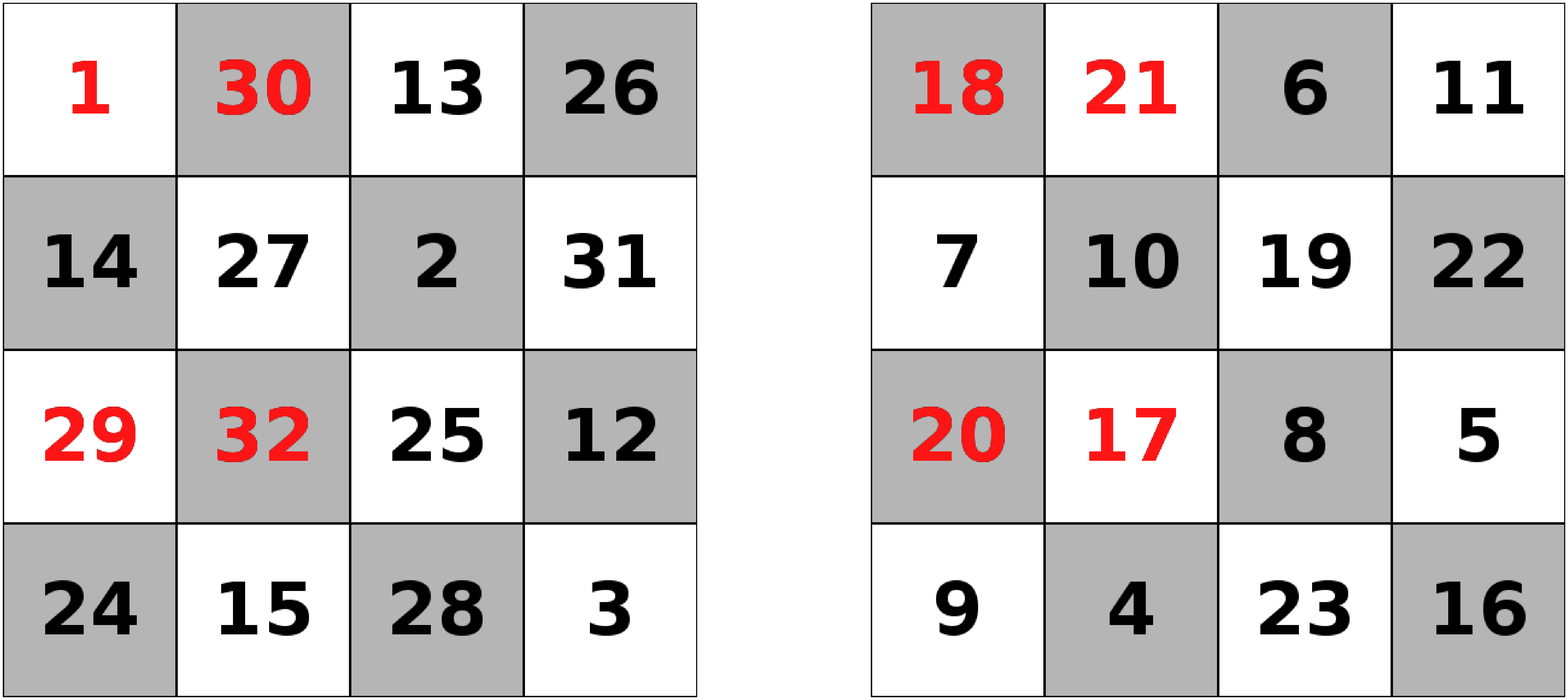}
\caption{A $4 \times 4 \times 2$ bi-sited tour}\label{f:4x4x2}
\end{figure}

\begin{figure}[ht]
\centering
\includegraphics[scale=0.2]{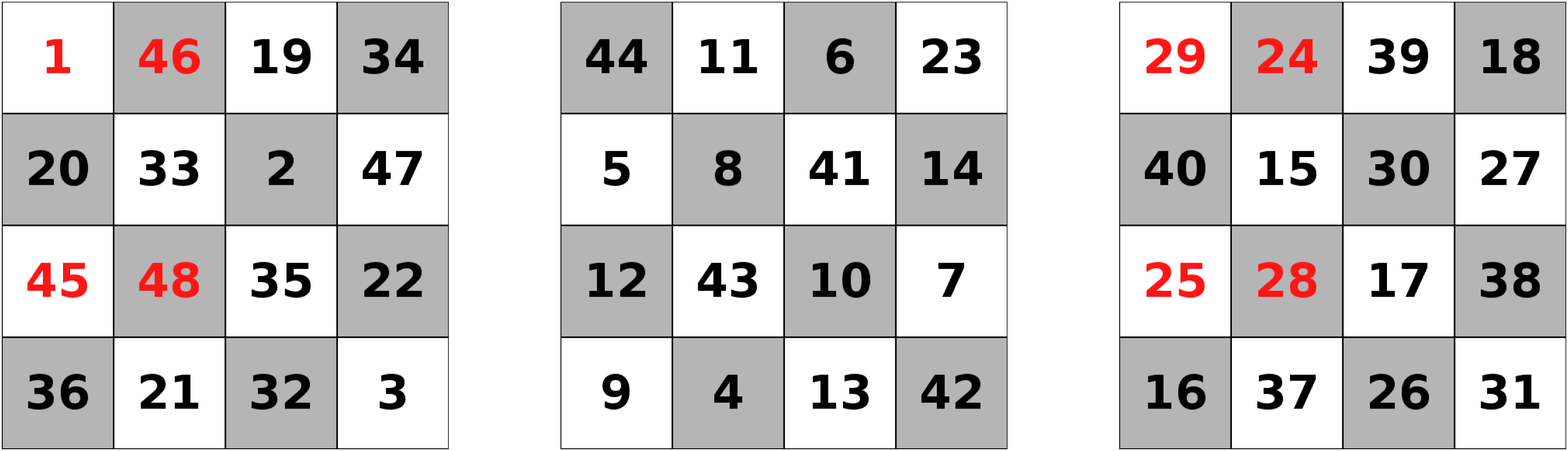}
\caption{A $4 \times 4 \times 3$ bi-sited tour by a knight move}\label{f:4x4x3}
\end{figure}

Then we can construct a closed tour of an $n \times m \times 2$ board
by putting two copies of the open tour on top of each other and adding
the lines $((n,m,1),(n,m-2,2))$ and $((n,m,2),(n,m-2,1))$. Then,
repeating the proof of Proposition \ref{p:ext}, we get bi-sited closed tours
for $(n+4p) \times (m+4q) \times 2$ boards. Finally, Proposition
\ref{p:gain} yields bi-sited closed tours for $(n+4p) \times (m+4q)
\times 2k$ boards for all $k\geq 1$. 
\\
\\
{\bf ii) Case $\bf 4 \times 4 \times k$, for $\bf k\geq 2$:}
We start by exhibiting a $4 \times 4 \times 2$ and a $4 \times 4
\times 3$ bi-sited tour in Figures \ref{f:4x4x2} and \ref{f:4x4x3}. 

Notice the sites in the top left corners of the top and bottom layers
of each of them, that is the lines $1-32$, $29-30$, $18-17$ and
$20-21$ in the $4 \times 4 \times 2$ tour and the lines 
$1-48$, $45-46$, $24-25$ and
$28-29$ in the $4 \times 4 \times 3$ tour. So we can
stack any number of these on top of each other to form bi-sited $4
\times 4 \times k$ tours for all $k$. More concretely we can form a $4
\times 4 \times 4$ tour by removing the line $20-21$ from a copy of
the $4 \times 4 \times 2$ tour and
placing it on top of another copy with the line $1-32$ removed, then
add in the lines $1-20$ and $32-21$. In a similar fashion we can add
any number of $4 \times 4 \times 2$ and $4 \times 4 \times 3$ tours
together.
\\
\\
{\bf iii) Case $\bf 4 \times 3 \times k$, for $\bf k\geq 2$:}
Again below we exhibit a $3
\times 4 \times 2$ and a $3 \times 4 \times 3$ bi-sittted tour. Recall
that they are equivalent to a $4 \times 3 \times 2$ and a $4 \times 3
\times 3$ tour. 
\begin{figure}[ht]
\centering
\includegraphics[scale=0.2]{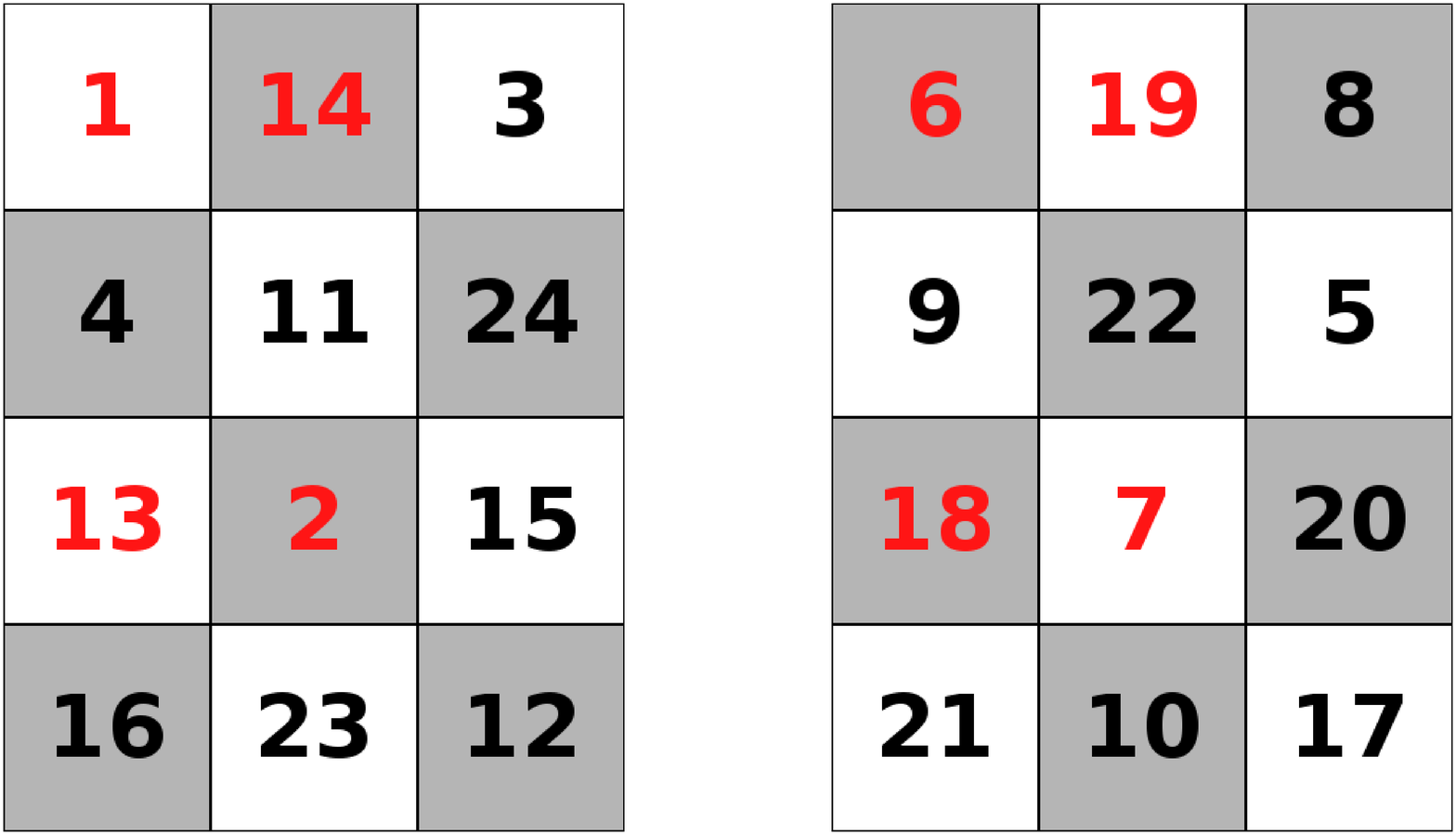}
\caption{A $4 \times 3 \times 2$ bi-sited tour}\label{f:4x3x2}
\end{figure}

\begin{figure}[ht]
\centering
\includegraphics[scale=0.2]{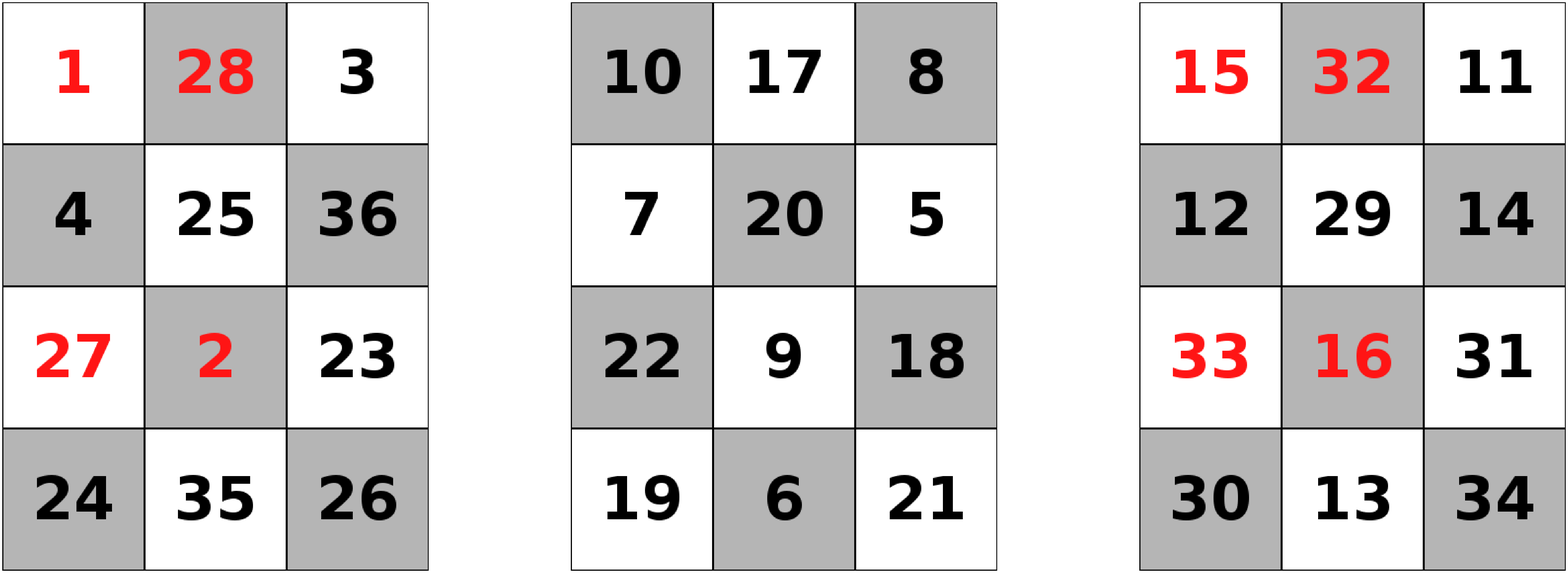}
\caption{A $4 \times 3 \times 3$ bi-sited tour}\label{f:4x3x3}
\end{figure}

Note that they have sites in the top left corners of the top and bottom
layers. Proceeding as in ii), we obtain bi-sited $4 \times 3
\times k$ tours for all $k \geq 2$. 
\\
\\
{\bf iii) Case $\bf 4 \times 2 \times k$, for $\bf k\geq 2$:}
As a $4 \times 2 \times 2$ tour does
not exist, we rely on a $4 \times 2 \times 3$, a $4 \times 2 \times 4$
and a $4 \times 2 \times 5$ bi-sited tour.

We construct a $4 \times 6 \times 2$ tour by stacking two copies
of the $4 \times 3 \times 2$ tour together, see Figure \ref{f:4x3x2},
and by  removing the $11-12$ line from the left copy and the $1-2$
line from the right copy and adding in the $11-11$ and $2-12$
lines. Inductively, we obtain bi-sited $4 \times k \times 2$ tours for
all $k \equiv 0 $ mod$(3)$. 

Similarly we place the $4 \times 3 \times 2$ tour from Figure
\ref{f:4x3x2} to the left of the 
$4 \times 4 \times 2$ tour from Figure \ref{f:4x4x2}. Then we 
remove the $1-32$ line 
from the $4 \times 4 \times 2$ tour and the $11-12$ from the $4 \times 3
\times 2$ tour and adding in the $11-1$ and $12-32$ lines. This settles
the case $k \equiv 1$ mod$(3)$.

\begin{figure}[ht]
\centering
\includegraphics[scale=0.2]{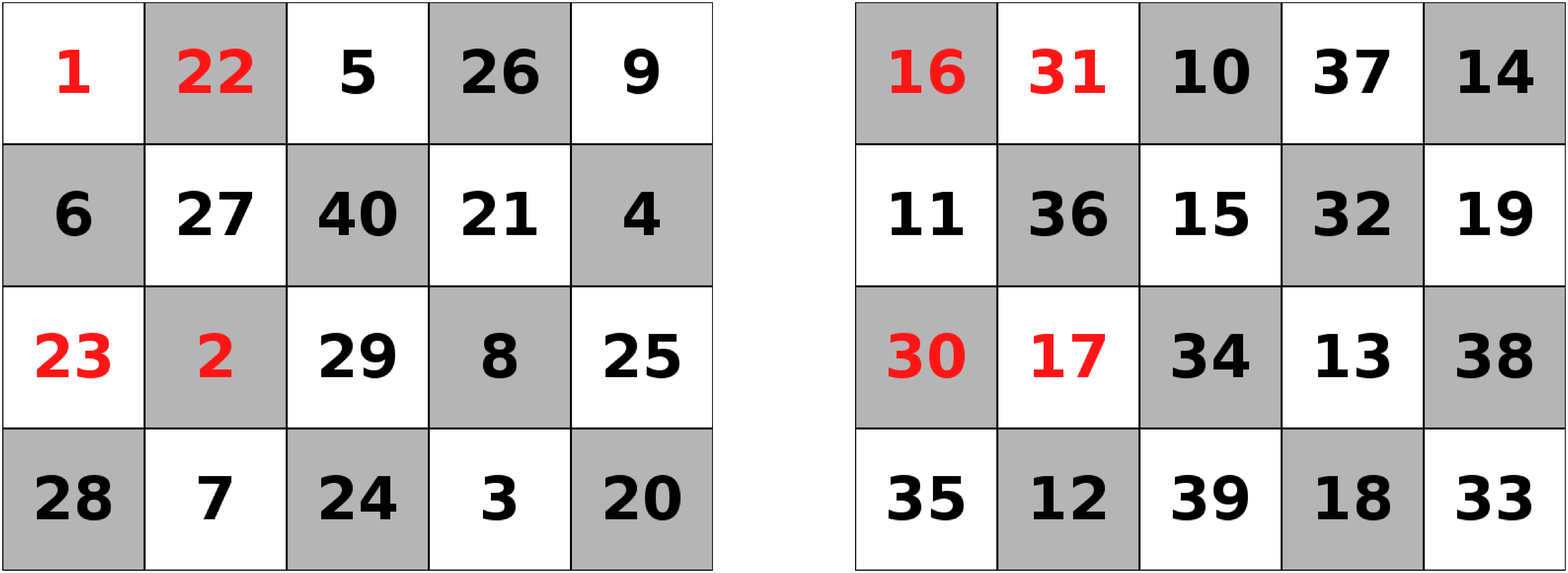}
\caption{A $4 \times 5 \times 2$ bi-sited tour}\label{f:4x5x2}
\end{figure}

Finally we give a bi-sited $5 \times 4 \times 2$ closed tour in Figure
\ref{f:4x5x2}. Repeating the procedure, this yields the case $k \equiv
2$ mod$(3)$ for $k\geq 5$. 

\noindent{\bf iv) Case $\bf 3 \times 2 \times k$, for $\bf k\geq 4$:} We now
construct a $3 \times 2 \times 8$ tour by stacking together two copies
of the $3 \times 2 \times 4$ tour, given in Figure \ref{f:4x3x2}, by
removing the line $15-16$ in the first copy and the line $8-9$ in the
second copy and adding in the lines $15-8$ and $16-9$. By induction we
derive $3 \times 2 \times k$ tours for  $k \equiv 0$
mod$(4)$. 

\begin{figure}[ht]
\centering
\includegraphics[scale=0.2]{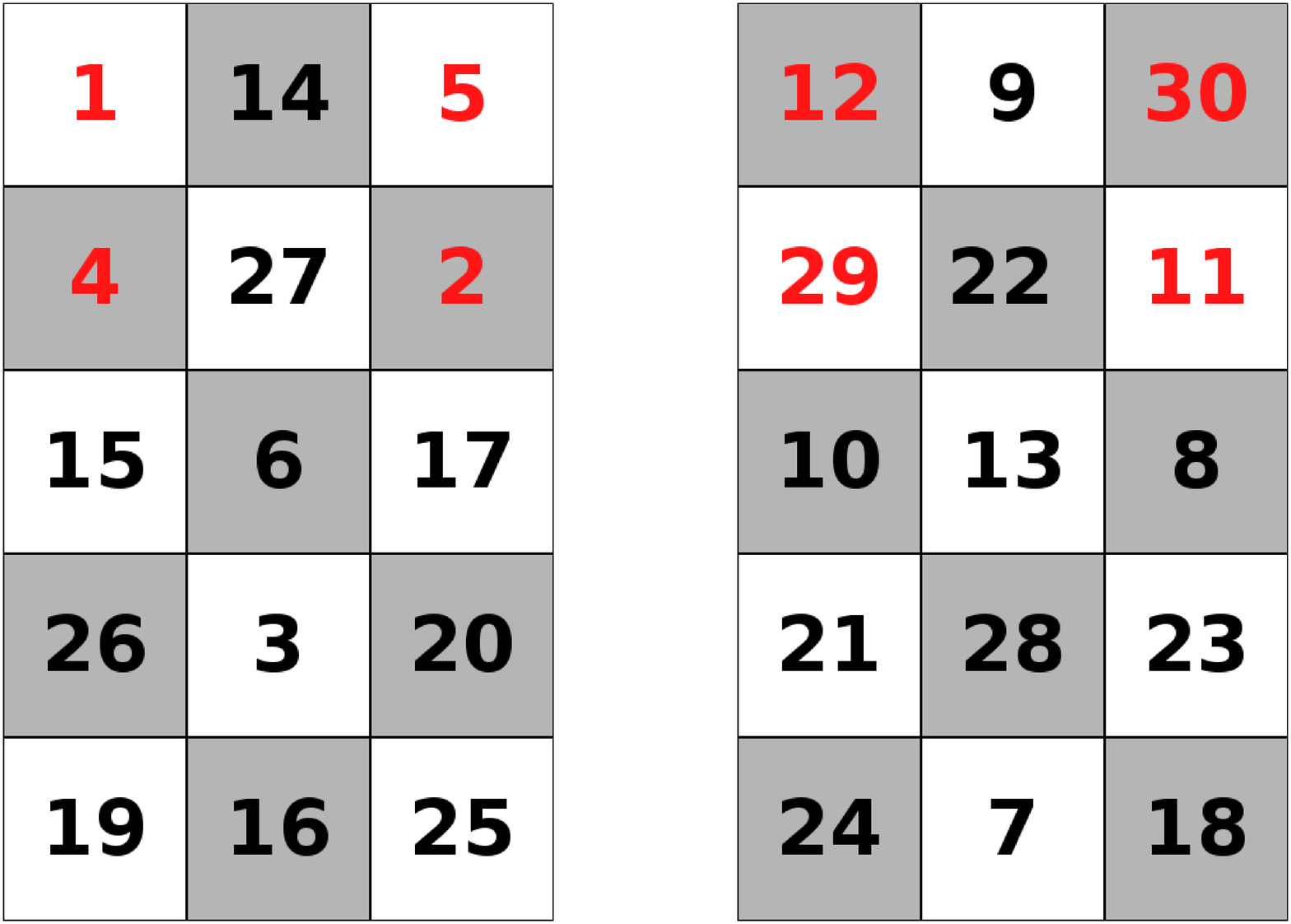}
\quad \includegraphics[scale=0.2]{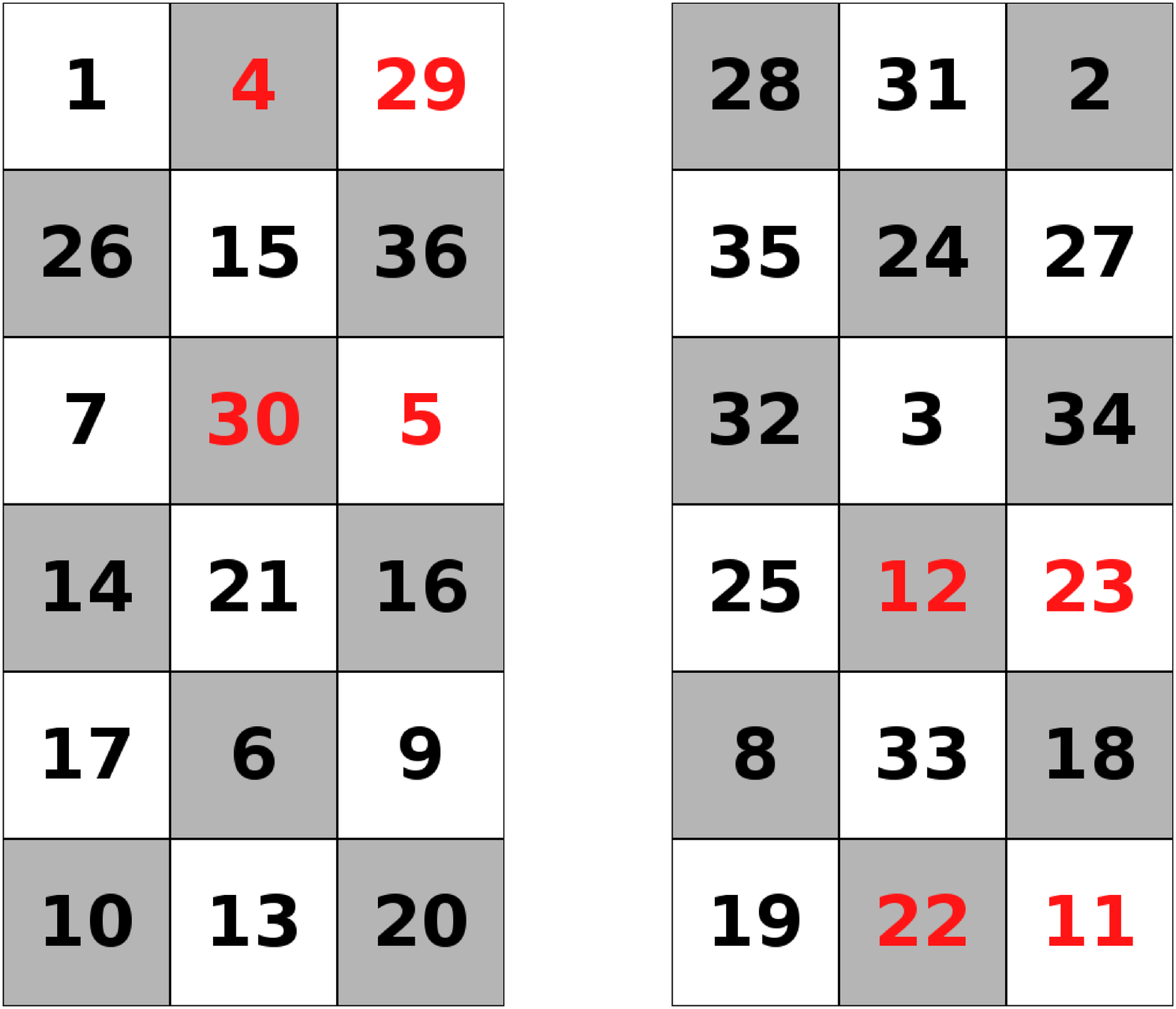} 
\caption{A $3 \times 5 \times 2$ and a $3 \times 6 \times 2$ bi-sited
  tour}\label{f:3x5x2} 
\end{figure}

\begin{figure}[ht]
\centering
\includegraphics[scale=0.2]{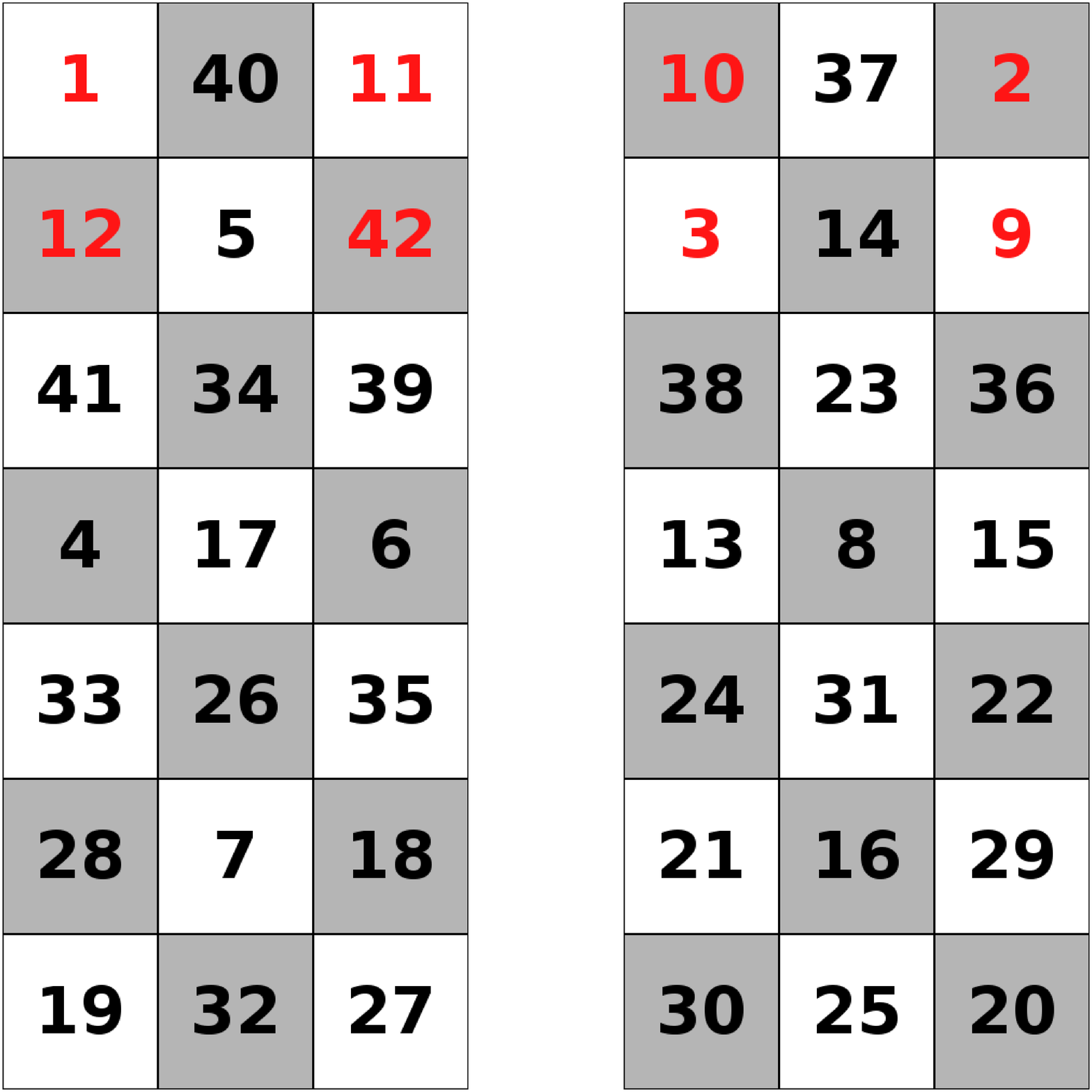}
\caption{A $3 \times 7 \times 2$ bi-sited tour}\label{f:3x7x2} 
\end{figure}

It remains exhibit tours of size $3 \times 5 \times
2$, $3 \times 6 \times 2$ in Figure \ref{f:3x5x2},  and $3 \times 7
\times 2$ in Figure \ref{f:3x7x2} , on which we 
consider the lines $(20-21)$ and $(23,24)$, $(28-29)$
and $(8,7)$, and $(33-34)$ and $(30,31)$ respectively.

\begin{figure}[ht]
\centering
\includegraphics[scale=0.2]{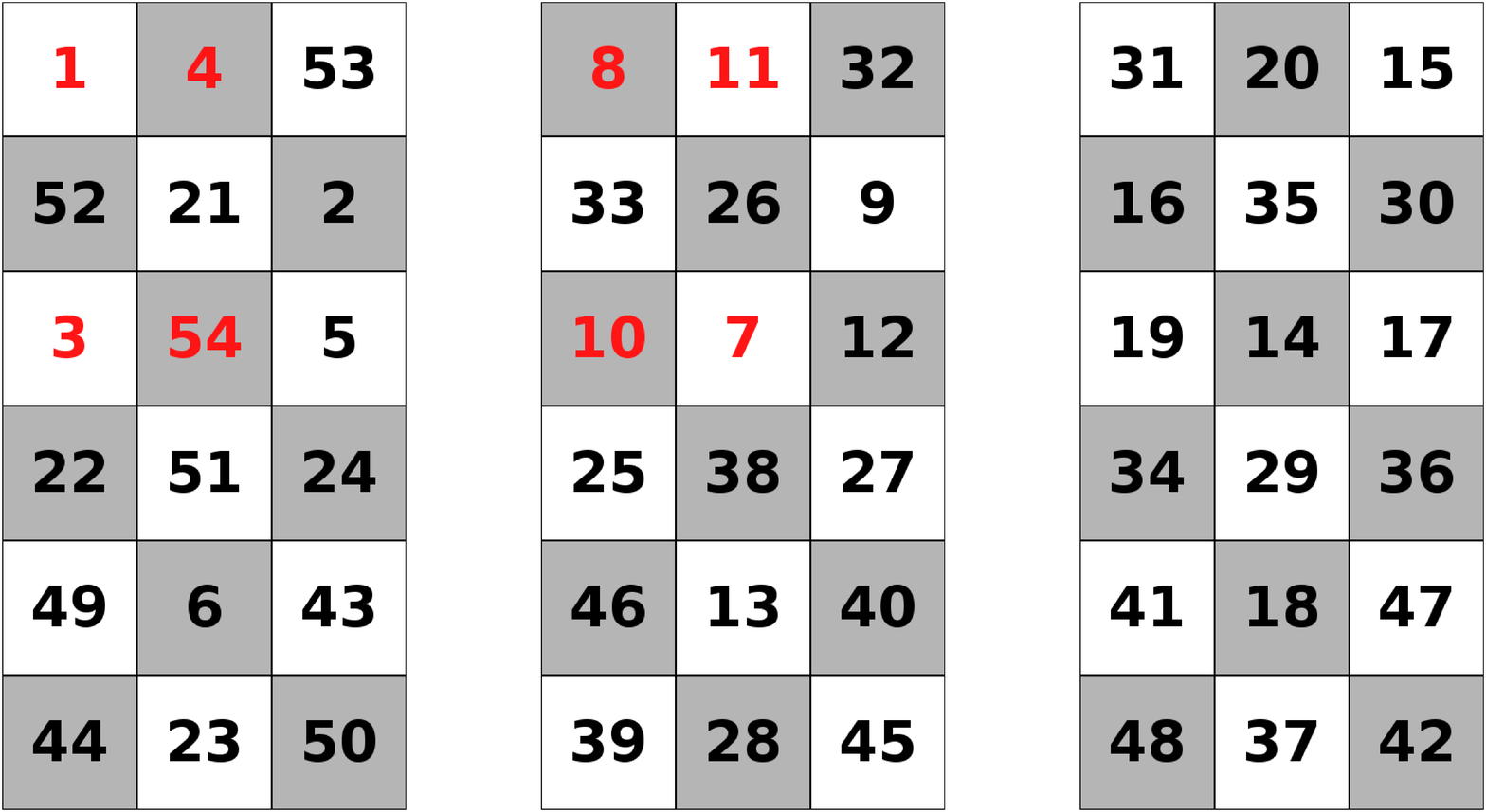}
\caption{A $3 \times 3 \times 6$ bi-sited tour.}\label{f:3x3x6} 
\end{figure}

\noindent{\bf v) Case $\bf 3 \times 3 \times 2k$, for $\bf k\geq 3$:}
For $k\geq 5$, it follows from Theorem \ref{t:S2} and Proposition
\ref{p:gain}. It remains $k=3$ and $k=4$. Firstly take 
the $4 \times 3 \times 3$ tour from Figure \ref{f:4x3x3} we can join
two of these together to form a $8 \times 3 \times 3$ 
tour by deleting the $23-24$ line in the first copy and the $7-8$ line
in the second copy and adding in the lines $7-24$ and $8-23$. We
conclude by giving a $3 \times 3 \times 6$ tour in Figure \ref{f:3x3x6}.

\newpage
\begin{thebibliography}{xxxxxx}
\bibitem[Cai]{Cai} G.\ Cairns: \emph{Pillow Chess}, Mathematics
  Magazine, 75, no.\ 3.\ (June, 2002), 173--186. 
\bibitem[De]{De} J. DeMaio, Which Chessboards have a Closed Knight's
  Tour within the Cube?, The Electron.\ J.\ Combin.\ 14,
  (2007).
\bibitem[DeM]{DeM} J.\ DeMaio and B.\  Mathew: \emph{Which chessboards
    have a closed knight's tour within the rectangular prism?} 
Electron.\ J.\ Combin.\ 18 (2011), no.\ 1, Paper 8, 14 pp.
\bibitem[Eul]{Eul} L.\ Euler; \emph{Solution d'une question curieuse
    qui ne paro\^it soumise \`a aucune analyse}, Histoire de
  l'Acad\'emie Royale des Sciences et des Belles-Lettres de Berlin, 
  vol.\ 15 (1759), 310--337.
\bibitem[Pab]{Pab} I.\ Parberry: \emph{An Efficient Algorithm for the Knight's Tour Problem}, Discrete Applied Mathematics, Vol.\ 73, pp.\ 251--260, 1997.
\bibitem[Kee]{Kee} M.R.\ Keen: \emph{The knight's tour},
  http://www.markkeen.com/knight/index.html 
\bibitem[Knu]{Knu} D.\ E.\ Knuth, \emph{Leaper Graphs}, Math.\ Gazette
  78 (1994), 274--297. 
\bibitem[Kum]{Kum} A.\ Kumar: \emph{Magic Knight's Tours in Higher
    Dimensions}, preprint arXiv:1201.0458.
\bibitem[QiW]{QiW} Y.\ Qing and J.J.\ Watkins: \emph{Knight's Tours for
    Cubes and Boxes}, Congressus Numerantium 181 (2006) 
    41--48.
\bibitem[Sch]{Sch} A.J.\ Schwenk: \emph{Which Rectangular Chessboards
    have a Knight's Tour?} Mathematics Magazine 64:5 (December 1991)
    325--332. 
\bibitem[Ste]{Ste} I.\ Stewart: \emph{Solid Knight's Tours}, Journal
  of Recreational Mathematics, Vol.\ 4 (1), January 1971.
\bibitem[Wat]{Wat} J.J.\ Watkins: \emph{Across the board: the
    mathematics of chessboard problems}, Princeton University Press,
    Princeton, NJ, (2004) xii+257 pp.\ ISBN: 0-691-11503-6. 
\bibitem[Wat1]{Wat1} J.J.\ Watkins, \emph{Knight's tours on cylinders
    and other surfaces}, Congr.\ Numer.\ 143 (2000), 117--127.
\bibitem[Wat2]{Wat2} J.J.\ Watkins, \emph{Knight's tours on a torus},
  Mathematics Magazine, 70:3 (1997), 175--184.
\bibitem[Wat3]{Wat3} J.J.\ Watkins, \emph{Across the board: the
    mathematics of chessboard problems}, Princeton University Press, 
    Princeton, NJ (2004). 
\bibitem[Fro]{Fro} M.\ Frolow: \emph{Les Carrés Magiques},(1886), Plate VII.
\end{thebibliography}
\end{document}